\documentclass{article}

\usepackage{arxiv}
\usepackage[utf8]{inputenc} 
\usepackage[T1]{fontenc}    
\usepackage{hyperref}       
\usepackage{url}            
\usepackage{booktabs}       
\usepackage{amsfonts}       
\usepackage{nicefrac}       
\usepackage{microtype}      
\usepackage{lipsum}
\usepackage{graphicx}
\graphicspath{ {./figs/} }

\usepackage{amsmath}
\usepackage{bm}
\usepackage{xcolor}
\usepackage{subfig}

\graphicspath{{./figs/}}
\usepackage{bm}
\def\Abold{\mathbf{A}}
 
\def\Bbold{\mathbf{B}} 
\def\bbold{\mathbf{b}} 
\def\Cbold{\mathbf{C}}

\def\dbold{\mathbf{d}}

\def\Fbold{\mathbf{F}} 
\def\fbold{\mathbf{f\/}} 
\def\Gbold{\mathbf{G}} 
 
\def\Hbold{\mathbf{H}}

\def\Kbold{\mathbf{K}} 
 
\def\Lbold{\mathbf{L}} 
 
\def\Mbold{\mathbf{M}} 
 
\def\Nbold{\mathbf{N}}

\def\pbold{\mathbf{p\/}}

\def\rbold{\mathbf{r}}

\def\tbold{\mathbf{t}} 
 
\def\ubold{\mathbf{u}}

\def\xbold{\mathbf{x}} 
 
\def\ybold{\mathbf{y}}

\def\barKbold{\mathbf{\bar K}}

\def\zerobold{\bm{0}}

\def\lambdabold{\bm{\lambda}}

\def\pibold{\bm{\pibold}}

\def\boxed#1{\setbox0=\hbox{$\displaystyle{#1}$}\hbox{\lower.6pt\hbox{\lower
 3pt\hbox{\lower 1\dp0\hbox{\vbox{\hrule height .6pt \hbox{\vrule width
 .6pt \hskip 3pt\vbox{\vskip 3pt\box0\vskip3pt}\hskip 3pt \vrule width
    .6pt}\hrule height .6pt}}}}}}
\def\dddot#1{{\mathop{#1}\limits^{\vbox to -1.4pt{\kern-2pt
   \hbox{\tenrm...}\vss}}}}
\def\ddddot#1{{\mathop{#1}\limits^{\vbox to -1.4pt{\kern-2pt
   \hbox{\tenrm....}\vss}}}}

\def\half{{\textstyle{1\over2}}}

\def\squeezedmatrixrt#1{\null\,\vcenter{\normalbaselines\m@th
    \ialign{\hfil$##$&&$\quad\!\!\!$\hfil$##$\crcr
      \mathstrut\crcr\noalign{\kern-\baselineskip}
      #1\crcr\mathstrut\crcr\noalign{\kern-\baselineskip}}}\,}
\def\crampedmatrix#1{\null\,\vcenter{\normalbaselines\m@th
    \ialign{\hfil$##$\hfil&&$\,\,$\hfil$##$\hfil\crcr
      \mathstrut\crcr\noalign{\kern-\baselineskip}
      #1\crcr\mathstrut\crcr\noalign{\kern-\baselineskip}}}\,}
\def\crampedmatrixrt#1{\null\,\vcenter{\normalbaselines\m@th
    \ialign{\hfil$##$&&$\,\,$\hfil$##$\crcr
      \mathstrut\crcr\noalign{\kern-\baselineskip}
      #1\crcr\mathstrut\crcr\noalign{\kern-\baselineskip}}}\,}

\def\bordertable#1{\begingroup \m@th
  \setbox\z@\vbox{\def\cr{\crcr\noalign{\kern2\p@\global\let\cr\endline}}%
    \ialign{$##$\hfil\kern2\p@\kern\p@renwd&\thinspace\hfil$##$\hfil
      &&\quad\hfil$##$\hfil\crcr
      \omit\strut\hfil\crcr\noalign{\kern-\baselineskip}%
      #1\crcr\omit\strut\cr}}%
  \setbox\tw@\vbox{\unvcopy\z@\global\setbox\@ne\lastbox}%
  \setbox\tw@\hbox{\unhbox\@ne\unskip\global\setbox\@ne\lastbox}%
  \setbox\tw@\hbox{$\kern\wd\@ne\kern-\p@renwd\left.\kern-\wd\@ne
    \global\setbox\@ne\vbox{\box\@ne\kern2\p@}%
    \vcenter{\kern-\ht\@ne\unvbox\z@\kern-\baselineskip}\,\right.$}%
  \null\;\vbox{\kern\ht\@ne\box\tw@}\endgroup}

\def\squeezedmorematrix#1{\null\,\vcenter{\normalbaselines\m@th
    \ialign{\hfil$##$\hfil&&$\quad\!\!\!\!$\hfil$##$\hfil\crcr
      \mathstrut\crcr\noalign{\kern-\baselineskip}
      #1\crcr\mathstrut\crcr\noalign{\kern-\baselineskip}}}\,}

\def\squeezedlessmatrix#1{\null\,\vcenter{\normalbaselines\m@th
    \ialign{\hfil$##$\hfil&&$\quad\!$\hfil$##$\hfil\crcr
      \mathstrut\crcr\noalign{\kern-\baselineskip}
      #1\crcr\mathstrut\crcr\noalign{\kern-\baselineskip}}}\,}

\def\btop#1{{\mathop{#1}\limits^{\vbox to -1.4pt{\kern-3pt
   \hbox{$\scriptscriptstyle\bullet$}\vss}}}}
\def\bbtop#1{{\mathop{#1}\limits^{\vbox to -1.4pt{\kern-3pt
   \hbox{$\scriptscriptstyle\bullet\bullet$}\vss}}}}
\catcode`@=12 

\chardef\@=`\@
\def\wiggle{\lower3pt\hbox{\char`\~}}

\def\diag{\hbox{\bf diag}}

\chardef\caret=`\^

\def\sqr#1#2{{\vcenter{\hrule height.#2pt
           \hbox{\vrule width.#2pt height#1pt \kern#1pt
           \vrule width.#2pt}
           \hrule height.#2pt}}}

\mathchardef\triangledown="0235

\mathchardef\varGamma="0100
\mathchardef\varDelta="0101
\mathchardef\varTheta="0102
\mathchardef\varLambda="0103
\mathchardef\varXi="0104
\mathchardef\varPi="0105
\mathchardef\varSigma="0106
\mathchardef\varUpsilon="0107
\mathchardef\varPhi="0108
\mathchardef\varPsi="0109
\mathchardef\varOmega="010A
\mathchardef\varkappa="017E

\def\overset#1\to#2{{\mathop{#2}^{#1}}}
\def\underset#1\to#2{{\mathop{#2}_{#1}}}
\def\oversetbrace#1\to#2{{\overbrace{#2}^{#1}}}
\def\undersetbrace#1\to#2{{\underbrace{#2}_{#1}}}

\def\@{\char'100 }

\def\Transp{^{\mathsf{T}}}

\title{Partitioned solution strategies for coupled BEM-FEM acoustic fluid-structure interaction problems

\thanks{\textit{\underline{Citation}}: 
\textbf{Luis Rodríguez-Tembleque, José A. González, Antonio Cerrato,
Partitioned solution strategies for coupled BEM–FEM acoustic fluid–structure interaction problems,
Computers \& Structures,
Volume 152,
2015,
Pages 45-58,
ISSN 0045-7949,
https://doi.org/10.1016/j.compstruc.2015.02.018.}} 
}

\author{
  Luis Rodríguez-Tembleque, José A. González, Antonio Cerrato  \\
  Escuela Ténica Superior de Ingeniería \\
  Universidad de Sevilla \\
  Camino de los Descubrimientos s/n, 41092 Sevilla, Spain\\
  \texttt{\{luisroteso, japerez, antoniocerrato\}@us.es} \\
}

\begin{document}
\maketitle
\begin{abstract}
Abstract: This paper investigates two FEM–BEM coupling formulations for acoustic fluid-structure interaction (FSI) problems, using the Finite Element Method (FEM) to model the structure and the Boundary Element Method (BEM) to represent a linear acoustic fluid. The coupling methods described interconnect fluid and structure using classical or localized Lagrange multipliers, allowing the connection of non-matching interfaces. First coupling technique is the well known mortar method, that uses classical multipliers and is compared with a new formulation of the method of localized Lagrange multipliers (LLM) for FSI applications with non-matching interfaces. The proposed non-overlapping domain decomposition technique uses a classical non-symmetrical acoustic BEM formulation for the fluid, although a symmetric Galerkin BEM formulation could be used as well. A comparison between the localized methodology and the mortar method in highly non conforming interface meshes is presented. Furthermore, the methodology proposes an iterative preconditioned and projected bi-conjugate gradient solver which presents very good scalability properties in the solution of this kind of problems.
\end{abstract}


\section{Introduction}\label{sec_1}
Reductions in noise emissions have high priority in the design
process of vibrating fluid-structure systems. These acoustic
fluid-structure interaction (FSI) problems are commonly found in
many engineering applications \cite{Ohayon95}, and the numerical
simulation of the interaction between the vibrating structure
gives fundamental information for optimizing the design of the
structure. In some situations one can perform the simulations
neglecting the influence of the acoustic field on the vibrating
structure. However, this is not acceptable for thin and flexible
structures that are easily excited by the acoustic pressure. For
these applications the acoustic field has to be fully coupled to
the vibrating structure.
%
The finite element method (FEM) have been applied to study this
kind of problems, many examples can be found in the book of Ohayon
and Soize \cite{Ohayon98} and Sandberg and Ohayon
\cite{Sandberg08a}.

The boundary element method (BEM) offers the major advantage over
the FEM, that only the boundary of the acoustic domain must be
discretized. Moreover, the Sommerfeld radiation condition for
exterior domains is inherently fulfilled, so it is especially more
appropriate than FEM to study exterior problems (ie. wave
propagation in infinite domains). For an introduction to the BEM,
it is referred to the monograph by Gaul et al. \cite{Gaul03}.
Coupling boundary element and finite element method in FSI, one
can benefit from the advantages of both numerical methodologies:
FEM is used to model the structure, and the BEM to model the
fluid. The first BEM-FEM coupling algorithm was developed by
Everstine and Henderson \cite{Everstine90}, and later, Chen et al.
\cite{Chen98} proposed a variational coupling scheme for Galerkin
methods. Further developments and applications of BEM-FEM methods
for structure-acoustic field interaction can be found in the works
of Gaul and Wenzel \cite{Gaul02}, Czygan and von Estorff
\cite{Czygan02}, and Langer and Antes \cite{Langer03a}, and  more
recently, Fritze et al. \cite{Fritze05}, Soares \cite{Soares09},
He et al. \cite{He12} and Soares and Godinho \cite{Soares12}.

%
Based in a mortar scheme, Fischer and Gaul \cite{Fischer05}
proposed an efficient FEM-BEM coupling in FSI which allows to
connect dissimilar meshes, using a to solve coupled
acoustical-fluid (BEM) structure (FEM), via classical Lagrange
multipliers. The mortar element method was originally introduced
by Bernardi et al. \cite{Bernardi94}. It offers the big advantage,
that non-conforming discretizations can be coupled. One obtains
great flexibility for meshing the subdomains and an increased
efficiency in the case that a coarser mesh is sufficient in one of
the subdomains.
%
Formulations, based on classical Lagrange multiplier fields, are
quite effective but tend to generate monolithic schemes that do
not preserve software modularity. To obtain a partitioned scheme,
Park and Felippa \cite{Park98,Park00,Park02} proposed a
formulation to connect non-matching FEM meshes. Non-matching
interfaces are treated by the method of localized Lagrange
multipliers (LLM), introducing a discrete surface \textit{frame}
interposed between the subdomains to approximate interface
displacements. This \textit{frame} is discretized and connected to
the BEM or FEM substructures by using LLM collocated at the
interface nodes. The application of BEM and FEM coupling in
elastostatics using localized Lagrange multipliers has been done
by Gonz\'{a}lez and Park \cite{Gonzalez07}, and the extension to
fluid-structure field interaction, by Park et al.
\cite{Park01,Park11}, Ross et al. \cite{Ross08, Ross09} and
Gonz\'{a}lez and Park \cite{Gonzalez12a} and Gonz\'{a}lez et al.
\cite{Gonzalez12b}.

In the mid-frequency regimes the acoustic fluid-structure problems
require fine meshes and, as a result, they generate a large number
of degrees of freedom. In this context, domain decomposition
methodologies (DDM) have appeared as a powerful numerical tool for
solving this large-scale systems. One of the most important
strategies is the finite element tearing and interconnecting
(FETI) method. The FETI methodology was proposed by Farhat and
Roux \cite{Farhat91} in the mid-90s, and it is an effective DDM
for the parallel solution of finite element problems partitioned
into subdomains. The global continuity across the subdomains
interfaces is enforced by classical Lagrange multipliers, which
leads to a saddle point problem that can be solved iteratively via
its dual problem. The dual problem leads to a linear flexibility
equations system for the Lagrange multipliers which is solved by a
preconditioned conjugate gradient (PCG) algorithm. The success of
the FETI method is due to its scalability with respect to the
problem size and number of subdomains \cite{Farhat00a,Farhat00b},
so the total solution time is approximately constant using a
smaller element size by multiplying proportionately the number of
processors. The numerical effort needed to solve the flexibility
system iteratively using a PCG algorithm is controlled by the
condition number of the system. Farhat et al.
\cite{Farhat00b,Farhat98a,Farhat98b} and Mandel and Tezaur
\cite{Mandel96} estimated the condition number of the system for
the FETI method as a function of the number of element per
subdomain: $\mathcal{O}((1+\log({L}/{h}))^2)$. Here, $L$ and $h$
denote the average size of the subdomains and finite elements,
respectively. Note that condition number is bounded independently
of the number of subdomains, so it is a necessary condition to
achieve numerical scalability. For some further development of the
FETI methodology (FETI-DP), the reader can turn to
\cite{Farhat00c,Farhat01}, and its application to acoustics FSI
problems, to Farhat et al. \cite{Farhat00d,Farhat00e} and Li et
al. \cite{Li11}.

The Boundary Element Tearing and Interconnecting (BETI) method
came up as a direct extension of the FETI to the BEM. The BETI
method was recently introduced by Langer and Steinbach
\cite{Langer03} as a counter part of the FETI methods. This
methodology extends the tearing and interconnecting technique to
symmetric Galerkin boundary element method (SGBEM)
\cite{Bonnet98,Gray98} in order to obtain symmetric system
matrices and, therefore, the use of FETI PCG solver becomes
feasible. The SGBEM is used to construct the
\emph{Steklov-Poincar\'e} operators instead of the finite element
based Schur complements, so the advantageous properties of FETI
methods remain valid for BETI methods as well. It can be
demonstrated numerically \cite{Langer03,Of09} that, using a
preconditioner with appropriated scaling matrices, the condition
number of the preconditioned BETI system is
$\mathcal{O}((1+\log({L}/{h}))^2)$, providing the same scalability
characteristics than the FETI method. This BETI methodology has
has been successfully applied to different kind of problems (i.e.
Bouchala et al. \cite{Bouchala09} proposed a BETI scheme for
contact problems), but it requires the implementation of SGBEM,
what is not straightforward. If it is replaced by a
non-symmetrical BEM formulation to approximate the
\emph{Steklov-Poincar\'{e}} operators of the floating
substructures, the flexibility equations become non-symmetric and
different solution strategies should be considered.

This work presents a non-symmetric FE-BETI formulation to solve
vibro-acoustic FSI problems based on
\cite{Park97,Filho97,Gonzalez12c}. The resulting non-symmetrical
flexibility system is solved by a new iterative solution procedure
based in a Bi-Conjugate Gradient Stabilized (Bi-CGSTAB) algorithm.
The formulation is presented together with some benchmark examples
that allow to compare the mortar and LLM schemes for
non-conforming approximations, and to check convergence and
accuracy of the non-symmetric FE-BETI algorithm as well as
scalability properties.

The paper is organised as follows. After the introduction, an
acoustics FSI partitioned formulation is developed. In Section 3,
special attention is paid to the choice of the coupling strategy:
mortar or LLM. The projected Bi-Conjugate Gradient Stabilized is
presented in Section 4, and scalability and convergence issues are
studied on different examples in Section 5. Finally, the paper
concludes with the summary and resulting conclusions.

\section{Acoustics FSI partitioned formulation}\label{sec_2}
A FEM structure and a BEM fluid domain are considered, so the
total virtual work of the system $\delta W_{T}$ can be expressed
as the addition of the virtual work done by the FEM structure
domain $\delta W_{S}$, the BEM fluid domain $\delta W_{F}$ and the
interface coupling contribution $\delta W_{C}$,

\begin{equation}
\delta W_{T} = \delta W_{S} + \delta W_{F} + \delta W_{C}
\end{equation}

\begin{figure*}
\centering
\subfloat[]{\includegraphics[width=5cm]{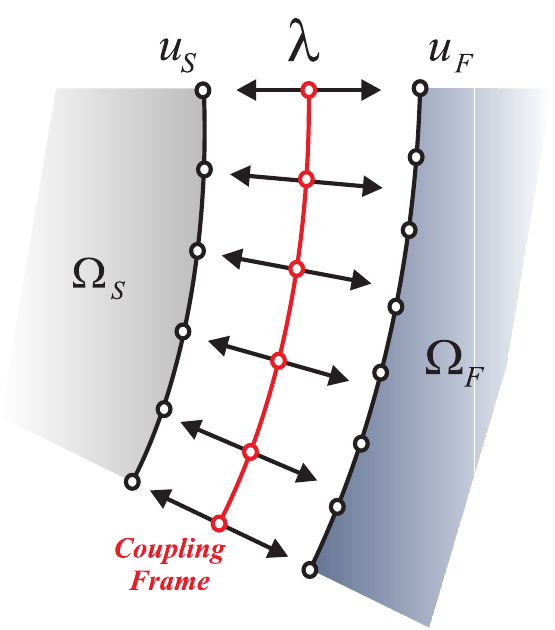}}
\hspace{2cm}
\subfloat[]{\includegraphics[width=5cm]{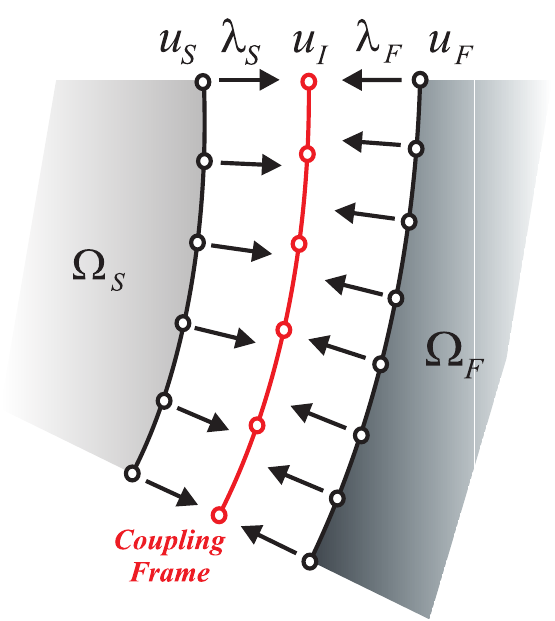}} 
\caption{Differences in the description of mortar (a) and localized (b) interfaces. Classical multipliers are used in (a) with a direct connection of the interfaces. In (b), an independent discretization of the interface is introduced and connected to the solid and fluid boundaries using localized Lagrange multipliers.}
\label{frame}
\end{figure*}

\subsection{Structural domain}
The virtual work of a flexible structure, $\delta W_{S}$, is
described by the principle of virtual work for a continuum body with
domain $\Omega_S$ and surface $\Gamma_S$ that, assuming small
displacements, can be written:
\begin{align}
\delta W_{S} =& \int_{\Omega_S} \sigma_{Sij} \delta u_{Si,j} d\Omega - \int_{\Omega_S} (\rho_S\omega^2 u_{Si} + b_{Si}) \delta u_{Si} d\Omega \nonumber \\
&- \int_{\Gamma_S} t_{Si}\delta u_{Si} d\Gamma
\end{align}
where $u_{Si}$ is the $i$-th component of the structural displacement
vector, a vector with the same number of components than the
dimension of the space; $\sigma_{Si,j}$ is the Cauchy stresses
tensor, $t_{Si}$ the applied surface tractions and $b_{Si}$ the
body forces. Finally, $\omega$ and $\rho_S$ are the angular
frequency of the harmonic oscillations and the density of the
structural material, respectively.

Next, the structure is discretized using the classical FEM
approximation, where the assembly of element contributions by the
direct stiffness method leads to the semidiscrete equation of
motion:
\begin{equation}\label{eq_PTV2}
\delta W_{S} = \delta \ubold_{S}\Transp \{ (\Kbold_S - \omega^2
\Mbold_S )\ubold_S - \fbold_S \}
\end{equation}
where $\Kbold_S$ is the stiffness matrix, $\Mbold_S$ is the mass
matrix, $\ubold_S$ is the vector of nodal displacements and
$\fbold_S$ are the applied nodal forces. Equation (\ref{eq_PTV2})
can be written in a more compact form:
\begin{equation}\label{VW_Structure}
\delta W_{S} = \delta \ubold_S\Transp \{\barKbold_S \ubold_S -
\fbold_S \}
\end{equation}
defining the dynamic stiffness matrix as $\barKbold_S =
\Kbold_S-\omega^2\Mbold_S$.

\subsection{Acoustic-fluid domain}
The governing equation for the linear a acoustic fluid with domain
$\Omega_F$ is the Helmholtz equation:
\begin{equation}\label{Helmholtz_eq}
\nabla^2 p_F + k_F^2 p_F = 0
\end{equation}
In this equation, $\nabla^2$ is the Laplacian, $p_F$ is the
acoustic fluid pressure, $k_F = \omega/c_F$ the wave number, $c_F$
the speed of sound in the fluid and $\omega$ the frequency.

On the fluid boundary $\Gamma_F = \Gamma_{Fu} \cup \Gamma_{Fr}$,
two types of boundary condition are considered:
\begin{itemize}
    \item [-] Neumann boundary condition:
            \begin{equation}\label{Neumann_BC}
                \frac{\partial p_F}{\partial n} = \rho_F \omega^2 u_{Fn} \quad \text{on} \quad \Gamma_{Fu}
            \end{equation}
    \item [-] Rigid boundary:
            \begin{equation}
                \frac{\partial p_F}{\partial n}=0 \quad \text{on} \quad \Gamma_{Fr}
            \end{equation}
\end{itemize}
where $n$ denotes the unit normal on the surface, $\rho_F$ is the
fluid density and $u_{Fn}$ represents the amplitude of the normal
displacement on the boundary.

The BEM formulation for a linear acoustic medium is well known and
can be found in many classical texts \cite{Wu00} and it is based
on the transformation of Helmholtz equation (\ref{Helmholtz_eq})
into a boundary integral equation. To do so, Helmholtz equation is
written in a weak form considering a weighted residual approach using the
Green's function $G(\xbold,\ybold)$ as the weighting function, being $\xbold$
the collocation point and $\ybold$ the source point. The
expression of the Green's function depends on the dimension of the
space, with:
\begin{equation}\label{G2D_eq}
G(\xbold,\ybold) = \frac{i}{4}H_{0}^{(1)}(k_F | \xbold-\ybold | )
\end{equation}
for two dimensions, where $H_{0}^{(1)}$ is the Hankel function of
the first kind and $i$ is the imaginary unit.

Applying Green's second theorem to the weighted residual and locating
the collocation point on the boundary, the resulting
boundary integral equation is:
\begin{align}
&C(\xbold)p_F(\xbold) + \int_{\Gamma_F} p_F(\ybold)\frac{\partial G(\xbold,\ybold)}{\partial n} d\Gamma = \nonumber \\
&\int_{\Gamma_F} G(\xbold,\ybold)\frac{\partial
p_F(\ybold)}{\partial n} d\Gamma \label{BIE_acoustics}
\end{align}
where $C(\xbold)$ is a coefficient that depends on the position of
point $\xbold$: $C(\xbold)=1$ for an internal point,
$C(\xbold)=\half$ for $\xbold$ on a smooth boundary $\Gamma_F$,
and $C(\xbold)=0$ for an external point.

Taking into account the Neumann boundary condition (\ref{Neumann_BC}), equation (\ref{BIE_acoustics}) can be rewritten as:
\begin{align}
&C(\xbold)p_F(\xbold) + \int_{\Gamma_F} p_F(\ybold)\frac{\partial G(\xbold,\ybold)}{\partial n} d\Gamma = \nonumber \\
&\rho_F \omega^2 \int_{\Gamma_{F_u}} G(\xbold,\ybold)
u_{Fn}(\xbold) d\Gamma \label{BIE_acoustics2}
\end{align}

Next, the BIE is discretized dividing the fluid boundary $\Gamma_F$ into
$n_{e}$ linear elements of surface $\Gamma_{e}$. Pressure and
displacement fields are approximated on each element $\Gamma_{e}$
by using linear shape functions:
\begin{equation}
p_F = \sum_{i=1}^{m} N_i p_{Fi} = \Nbold \pbold_{F}
\label{p_discrete}
\end{equation}
\begin{equation}
u_{Fn}= \sum_{i=1}^{m} N_i u_{Fni} = \Nbold \ubold_{F}
\label{u_discrete}
\end{equation}
where $p_{Fi}$ and $u_{{Fn}i}$ are the nodal values of acoustic
pressure and fluid normal displacement at node $i$, and $\Nbold$
is the shape function approximation matrix. A discrete boundary
integral equation is then obtained substituting equation
(\ref{p_discrete}) into equation (\ref{BIE_acoustics2}) and considering
that the point $\xbold$ is collocated on a boundary node:
\begin{align}
&C_i \delta_{ij}p_{Fj} + \sum_{e=1}^{n_{e}} \int_{\Gamma_e} \frac{\partial G(\xbold,\ybold)}{\partial n} N_j p_{Fj} d\Gamma_e = \nonumber \\
&\rho_F \omega^2 \sum_{e=1}^{n_{e}}\int_{\Gamma_e}
G(\xbold,\ybold)N_j u_{Fnj} d\Gamma_e \label{BIE_discrete}
\end{align}
being $\delta_{ij}$ is the Kronecker $\delta$-function. Equation
(\ref{BIE_discrete}) can be written in matrix form as:
\begin{equation}
\Hbold \pbold_{F} = \Gbold \ubold_{F} \label{BEM_eq}
\end{equation}
with the following definition for the matrix components:
\begin{equation}
H_{ij} = C_i \delta_{ij} + \sum_{e=1}^{n_{e}} \int_{\Gamma_e}
\frac{\partial G(\xbold,\ybold)}{\partial n} N_j d\Gamma_e
\end{equation}
\begin{equation}
G_{ij} = \rho_F \omega^2 \sum_{e=1}^{n_{e}}\int_{\Gamma_e}
G(\xbold,\ybold)N_j d\Gamma_e
\end{equation}

The virtual work of a BEM fluid subdomain can then be computed using a
weak statement for dynamic equilibrium reduced to the boundary.
This is done with Clapeyron equation \cite{Bonnet04,Gonzalez07}
that is expressed in the following form:
\begin{equation}\label{eqn_deltaWbem}
\delta W_{F} = \int_{\Gamma_F} ( p_{F} - t_{F}) \delta u_{Fn}
d\Gamma
\end{equation}
defining the external normal tractions imposed on the boundary as
$t_{F}$ and where the fluid pressure $p_F$ satisfies equation
\eqref{BIE_acoustics}.

Discretizing equation \eqref{eqn_deltaWbem} using the same BEM
mesh utilized for the fluid, a discrete approximation of the
virtual work is obtained:
\begin{align}
\delta W_{F} &= \delta \ubold_{F}\Transp  [ \int_{\Gamma_F} \Nbold\Transp\Nbold \, d\Gamma ] \{ \pbold_{F} - \tbold_{F} \} \nonumber \\
&= \delta \ubold_{F}\Transp \Mbold \{ \pbold_{F} - \tbold_{F} \}
\label{eqn_deltaWFdisc}
\end{align}
with a lumping matrix
\begin{equation}
\Mbold = \int_{\Gamma_F} \Nbold\Transp\Nbold \, d\Gamma
\end{equation}
that transforms distributed tractions into equivalent nodal
forces.

Substituting the discrete fluid pressures $\pbold_F$ coming from the BE equation
\eqref{BEM_eq} into the variational form \eqref{eqn_deltaWFdisc}, a final expression for the discrete variation is
obtained:
\begin{equation}
\label{eq_PTVF} \delta W_{F} = \delta \ubold_{F}\Transp\Mbold \{
\rho_F \omega^2 \Hbold^{-1}\Gbold\ubold_{F} - \tbold_{F} \}
\end{equation}
Note that this variational statement, obtained from a boundary
integral formulation, in general does not derive from an energy
functional and will be non-symmetric.

By comparison with equation \eqref{VW_Structure}, we conclude that the
equations for the fluid and the structure can be written using the
same expression:
\begin{equation}\label{VW_Fluid}
\delta W_{F} = \delta \ubold_{F}\Transp \{ \barKbold_{F}
\ubold_{F} - \fbold_{F} \}
\end{equation}
simply by defining an equivalent dynamic stiffness matrix for the
fluid $\barKbold_{F} = \rho_{F}\omega^2 \Mbold\Hbold^{-1}\Gbold$
and the vector of given external forces as $\fbold_{F} =
\Mbold\tbold_{F}$.

\section{Coupling strategies}\label{sec_3}
Two different strategies are investigated in this Section for the connection
of the fluid and the structure: \emph{Mortar} method and the
method of \emph{localized Lagrange multipliers}. A general
formulation is derived first for both methodologies and then they are compared
using a classical FSI example from \cite{Sandberg08b}.

\subsection{Mortar method}
We consider two coupled subdomains, $\Omega_S$ and $\Omega_F$,
sharing a common interface $\Gamma_C$. In Mortar methods, the work
associated to the tying interface will enforce the coupling
condition in a weak sense through the following expression:
\begin{equation}
W_{C} = \int_{\Gamma _{C}} \lambda (u_{Sn} - u_{Fn}) \, d\Gamma
\end{equation}
where $\lambda$ is the Lagrange Multiplier traction field over the
coupling interface $\Gamma_{C}$, and $u_{Sn}$ and $u_{Fn}$ are the
structure and fluid normal displacement fields over $\Gamma_{C}$.
The variation of this form leads to:
\begin{align}
&\delta W_{C} = \int_{\Gamma _{C}} \delta\lambda (u_{Sn}-u_{Fn}) \ d\Gamma + \int_{\Gamma_{C}} \delta u_{Sn} \lambda \ d\Gamma \nonumber \\
&- \int_{\Gamma _{C}} \delta u_{Fn} \lambda \ d\Gamma
\label{eq_Wc_Mortar}
\end{align}

The interface normal tractions $\lambda$ and the boundary normal
displacements $u_n$ of each domain ($S$, $F$) are interpolated on the boundary as
follows:
\begin{equation}
\lambda = \sum_{i=1}^{n_{I}} \hat{N}_{i} \lambda _{i}
\label{eq_discretiLambda}
\end{equation}
\begin{equation}
u_{Sn} = \sum_{i=1}^{n_{S}} N_{Si} u_{Sni}, \hspace{0.2cm} u_{Fn}
= \sum_{i=1}^{n_{F}} N_{Fi} u_{Fni} \label{eq_discretiDisp}
\end{equation}
where the shape functions $(N_{Si},N_{Fi})$ are defined
independently for the solid and the fluid side,
$(u_{Sni},u_{Fni})$ are the structure and fluid normal nodal
displacements and $(n_{I},n_{F},n_{S})$ are the
number of interface, fluid and structure boundary nodes.

\begin{figure}
\begin{center}
\includegraphics[width=12cm]{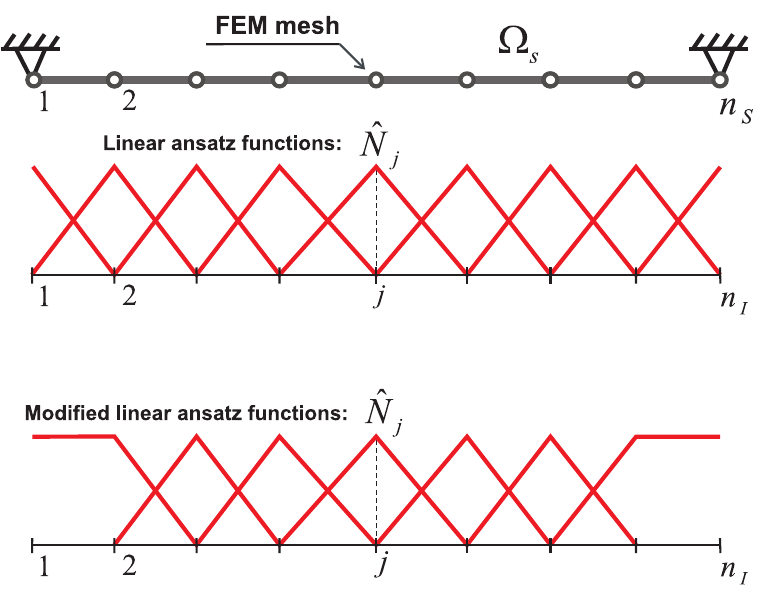}\\
\caption{Modified linear ansatz functions used with Mortar method
in the presence of Dirichlet boundary conditions for the
approximation of interface Lagrange
multipliers.}\label{Fig_ansatz}
\end{center}
\end{figure}

Lagrange multipliers are approximated using linear shape functions
$\hat{N}_{i}$ with the support of the discretisation on the solid
non-mortar side, same approximation than in \cite{Bernardi94}.
When Dirichlet boundary conditions exist on the boundary of
$\Gamma_{C}$, the shape functions $\hat{N}_{i}$ have to be
modified to avoid over-constrained conditions at those
restricted edges \cite{Wohlmuth00}, as represented in Figure
\ref{Fig_ansatz}.

\begin{figure*}
\centering
\includegraphics[width=6cm]{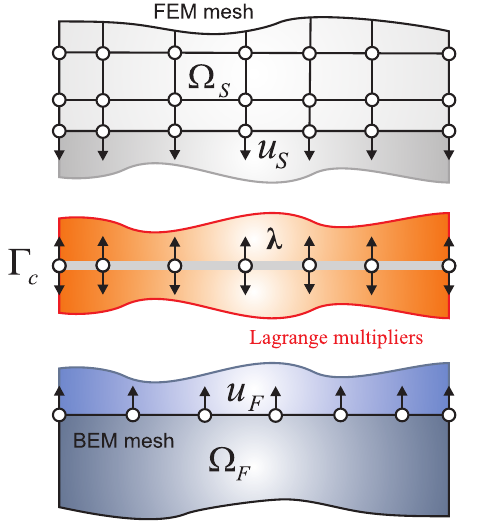} \includegraphics[width=7cm]{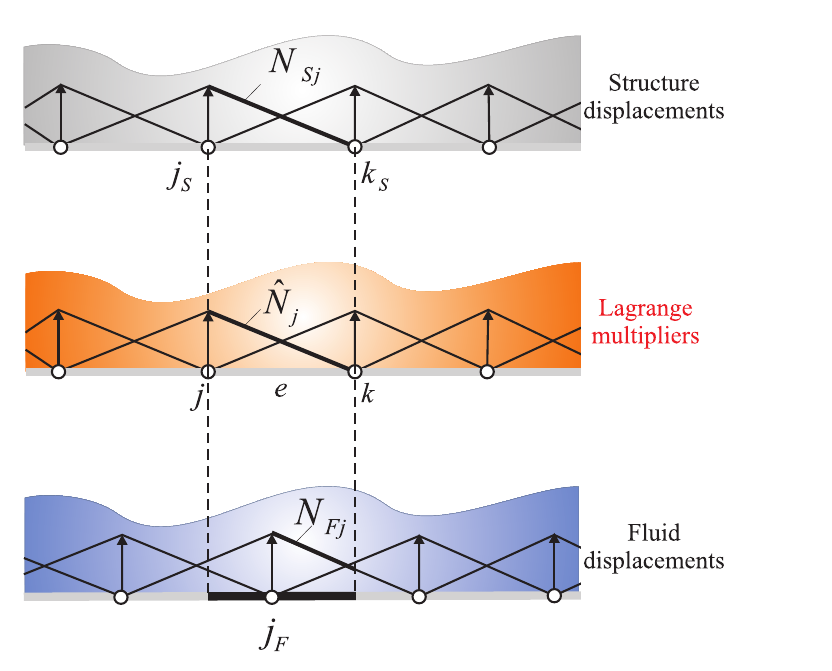}
\caption{Coupling with the Mortar method fluid and solid
interfaces. Approximation spaces for boundary displacements and
multipliers.}\label{Fig_ansatz_int}
\end{figure*}

Substituting approximations \eqref{eq_discretiLambda} and
\eqref{eq_discretiDisp} in the mortar coupling equation
\eqref{eq_Wc_Mortar}, the boundary integrals can be approximated:
\begin{align}
&\int_{\Gamma_{C}}\delta u_{Sn} \lambda\ d\Gamma = \nonumber \\
&\sum\limits_{e=1}^{n_{e}} \sum\limits_{i=1}^{n_{F}}
\sum\limits_{j=1}^{n_{S}} \delta u_{Sn i} [
\int_{\Gamma_{e}}N_{Si} \hat{N}_{j} \, d\Gamma ] \lambda_{j}
\label{mortar_1_expres}
\end{align}
\begin{align}
&\int_{\Gamma_{C}}\delta u_{Fn} \lambda\ d\Gamma =  \nonumber \\
&\sum\limits_{e=1}^{n_{e}} \sum\limits_{i=1}^{n_{F}} \sum\limits_{j=1}^{n_{S}}
\delta u_{Fn i} [ \int_{\Gamma_{e}}N_{Fi} \hat{N}_{j} \, d\Gamma
]\lambda_{j} \label{mortar_2_expres}
\end{align}
where $n_{e}$ is the number of elements on the non-mortar side and
$\Gamma _{e}$ is the element $e$ boundary (see
Figure \ref{Fig_ansatz_int}). The boundary integrals of equations
\eqref{mortar_1_expres} and \eqref{mortar_2_expres} are
assembled into matrices $\Abold_{S}$ and $\Abold_{F}$, so
equation \eqref{eq_Wc_Mortar} can be written in matrix form as follows:
\begin{align}
\delta W_{C} =& \delta\lambdabold\Transp (\Abold_S\Transp \ubold_{Sn} + \Abold_F\Transp \ubold_{Fn}) \nonumber \\
& + (\delta\ubold_{Sn})\Transp \Abold_S \lambdabold +
(\delta\ubold_{Fn})\Transp \Abold_F \lambdabold
\label{eq_coupling_mortar}
\end{align}
being $\lambdabold$, $\ubold_{Sn}$ and $\ubold_{Fn}$  the vectors of nodal
tractions , and nodal structure and fluid interface
normal displacements.


Normal displacement vectors on the boundary, $\ubold_{Sn}$ and
$\ubold_{Fn}$, can be obtained from the global vectors of structure
and fluid displacements:
\begin{equation}\label{displacement_local}
\ubold_{Sn} = \Bbold_{S}\Transp\ubold_{S}, \hspace{0.5cm}
\ubold_{Fn} = \Bbold_{F}\Transp\ubold_{F}
\end{equation}
where $\Bbold_{S}$ and $\Bbold_{F}$ are boolean matrices.
Substituting \eqref{displacement_local} into
\eqref{eq_coupling_mortar} and defining the coupling matrices:
\begin{equation}
\Cbold_{S} = \Bbold_{S} \Abold_{S}, \hspace{0.5cm} \Cbold_{F} =
\Bbold_{F} \Abold_{F}
\end{equation}
we arrive to the following expression for the virtual work of the interface:
\begin{equation}\label{VW_coupling_mortar}
\delta W_{C}= \delta\lambdabold\Transp (\Cbold_{S}^T \ubold_{S} +
\Cbold_{F}\Transp\ubold_{F}) + \delta\ubold_{S}\Transp \Cbold_{S}
\lambdabold + \delta\ubold_{F}\Transp \Cbold_{F} \lambdabold
\end{equation}

Total virtual work of the system $\delta W_{T}$ is then derived
from \eqref{VW_Structure}, \eqref{VW_Fluid} and
\eqref{VW_coupling_mortar} as:
\begin{align}
&\delta W_{T} = \delta \ubold_{S}\Transp\{\barKbold_{S}\ubold_{S} + \Cbold_{S} \lambdabold - \fbold_{S}\} + \nonumber \\
&\delta \ubold_{F}\Transp\{\barKbold_{F}\ubold_{F} + \Cbold_{F}
\lambdabold - \fbold_{F}\} + \delta
\lambdabold\Transp\{\Cbold_{S}\Transp \ubold_{S} +
\Cbold_{F}\Transp  \ubold_{F} \}
\end{align}
and from the stationary-point condition of this virtual work, the
following partitioned FSI equation set is obtained:
\begin{equation}\label{eqn_Mortarsystem}
    \begin{bmatrix}
    \barKbold_{S} & \zerobold                    & \Cbold_{S} \\
    \zerobold            & \barKbold_{F}          & \Cbold_{f} \\
    \Cbold_{S}\Transp       & \Cbold_{F}\Transp              & \zerobold
    \end{bmatrix}
    \begin{bmatrix}
    \ubold_{S} \\
    \ubold_{F}\\
    \lambdabold \\
    \end{bmatrix}
    =
    \begin{bmatrix}
    \fbold_{S} \\
    \fbold_{F} \\
    \zerobold
    \end{bmatrix}
\end{equation}

In general, for $p=1\dots n_{p}$ fluid and structure partitions,
equation \eqref{eqn_Mortarsystem} can be expressed in condensed
form as:
\begin{equation}\label{Mortar_system}
\begin{bmatrix}
\Kbold   & \Cbold \\
\Cbold\Transp& \zerobold
\end{bmatrix}
\begin{bmatrix}
\ubold \\
\lambdabold
\end{bmatrix}
=
\begin{bmatrix}
\fbold \\
\zerobold
\end{bmatrix}
\end{equation}
by simply defining the block-matrices:
\begin{equation}
\Kbold = \diag \begin{bmatrix} \Kbold_{\star}^{(1)}
\dots\Kbold_{\star}^{(n_{p})} \end{bmatrix}, \quad \Cbold =
\begin{bmatrix} \Cbold_{\star}^{(1)} \\ \vdots \\
\Cbold_{\star}^{(n_{p})} \end{bmatrix}
\end{equation}
and block-vectors:
\begin{equation}
\ubold = \begin{Bmatrix} \ubold_{\star}^{(1)} \\ \vdots \\
\ubold_{\star}^{(n_{p})} \end{Bmatrix}, \quad \fbold =
\begin{Bmatrix} \fbold_{\star}^{(1)} \\ \vdots \\
\fbold_{\star}^{(n_{p})} \end{Bmatrix}
\end{equation}
with subscript $\star = S,F$ indicating the type of model
associated to substructure $p$, i.e., (S) for solid modeled using
FEM or (F) for acoustic fluid using BEM.

After this reorganization and eliminating the displacements $\ubold$
from the first row of (\ref{Mortar_system}) using the relation:
\begin{equation}
\ubold = \Kbold^{-1}(\fbold - \Cbold \lambdabold)
\end{equation}
a compact non-symmetrical flexibility system is obtained:
\begin{equation}
\Fbold_{bb}\lambdabold = \bbold
\end{equation}
being $\Fbold_{bb} = \Cbold\Transp\Kbold^{-1}\Cbold$ a boundary
flexibility matrix and $\bbold=\Cbold^T\Kbold^{-1}\fbold$ the free
term.

\subsection{Localized Lagrange multipliers method}
The virtual work for the interface frame $\delta W_{C}$ can be
also evaluated applying the variationaly-based formulation proposed
by Park and Felippa \cite{Park98,Park00} and Gonz\'alez et al.
\cite{Gonzalez07}. The virtual work of the total system
$\delta W_{T}$ consists of contributions from the FE structure and BE
fluid, $\delta W_{S}$ and $\delta W_{F}$, plus the interface
frame $\delta W_{C}$. This formulation enforces the kinematical
positioning of the frame in a weak sense with the following
expression:
\begin{equation}\label{eq_int_LLM}
W_{C}=\int_{\Gamma _{C}} \{ \lambda_{S} (u_{Sn} - u_{In}) + \{
\lambda_{F} (u_{Fn}-u_{In})\}d\Gamma
\end{equation}%
where both integrals are extended to the boundary interface
$\Gamma_{C}$. The localized Lagrange multipliers and the
displacements of the structure interface are represented by
($\lambda_{S}$, $u_{Sn}$), and the fluid localized Lagrange
multipliers and displacements by ($\lambda_{F},u_{Fn}$). Finally,
the frame displacements are represented by $u_{In}$.

\begin{figure}
\begin{center}
  \includegraphics[width=7.3cm]{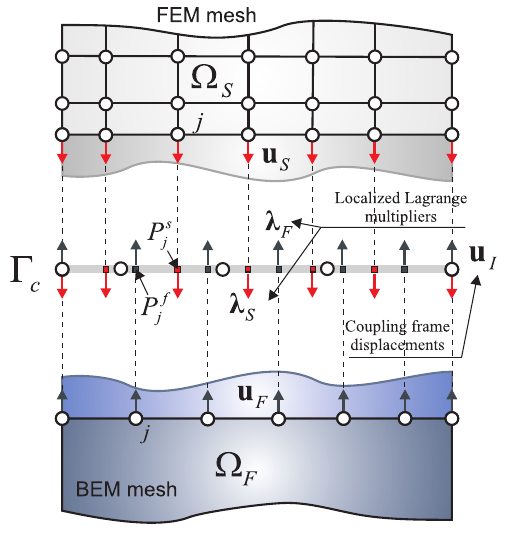}
  \caption{FSI BEM-FEM system with intercalated frame and localized Lagrange multipliers.}\label{Fig_ansatz_LLM}
\end{center}
\end{figure}

Equation \eqref{eq_int_LLM} can be written in matrix form as:
\begin{equation}
\delta W_{C} = \delta\{\lambdabold_S\Transp (\Bbold_{S}\Transp
\ubold_{S}-\Lbold_{S}\ubold_{I})\} + \delta\{\lambdabold_F\Transp
(\Bbold_{F}\Transp \ubold_{F}-\Lbold_{F}\ubold_{I})\}
\end{equation}
using two linear operators, $\Bbold_{S}$ to extract the structural
boundary displacements projected into the normal direction
and $\Bbold_{F}$ to extract the fluid boundary displacements.
$\Lbold_{S}$ and $\Lbold_{F}$ \cite{Park02,Gonzalez07} are matrices whose terms are obtained by evaluating the frame shape functions at the interface nodal position of the structure
and fluid $P_j^s$ and $P_j^f$ (see Figure \ref{Fig_ansatz_LLM}).

The total virtual work of the coupled BEM-FEM-Frame system can finally
be expressed as:
\begin{align}
\delta W_{T} &= \delta \ubold_{S}\Transp\{\barKbold_{S}\ubold_{S} + \Bbold_{S} \lambdabold_S - \fbold_{S}\} + \nonumber \\
&\delta \ubold_{F}\Transp\{\barKbold_{F}\ubold_{F} + \Bbold_{F} \lambdabold_F - \fbold_{F}\} + \delta \lambdabold_S\Transp\{\Bbold_{S}\Transp \ubold_{S} - \Lbold_S\ubold_{I} \} \nonumber \\
& + \delta \lambdabold_F\Transp\{\Bbold_{F}\Transp \ubold_{F} -
\Lbold_F\ubold_{I} \} + \delta \ubold_{I}\Transp \{\Lbold_S\Transp
\lambdabold_S+\Lbold_F\Transp \lambdabold_F \}
\end{align}

The stationarity condition of this variational form provides the
equations of motion, defined by the following system:
\begin{equation}\label{eqn_LLMsystem}
\begin{bmatrix}
\barKbold_{S}            & \zerobold                    & \Bbold_{S}    & \zerobold         & \zerobold \\
\zerobold                  & \barKbold_{F}             & \zerobold      & \Bbold_{F}       & \zerobold  \\
\Bbold_{S}\Transp             & \zerobold                    & \zerobold      & \zerobold         & \Lbold_S \\
\zerobold                  & \Bbold_{F}\Transp               & \zerobold      & \zerobold         & \Lbold_F \\
\zerobold                  & \zerobold                    &
\Lbold_S\Transp & \Lbold_F\Transp    & \zerobold
\end{bmatrix}
\begin{bmatrix}
\ubold_{S} \\
\ubold_{F}\\
\lambdabold_S \\
\lambdabold_F \\
\ubold_{I}\\
\end{bmatrix}
=
\begin{bmatrix}
\fbold_{S} \\
\fbold_{F} \\
\zerobold \\
\zerobold \\
\zerobold
\end{bmatrix}
\end{equation}

In general, if we have $n_{p}$ different fluid and structure
partitions, equation \eqref{eqn_LLMsystem} can be written in a
more compact form:
\begin{equation}\label{LLM_system}
\begin{bmatrix}
\Kbold   & \Bbold   & \zerobold \\
\Bbold\Transp& \zerobold   & \Lbold \\
\zerobold   & \Lbold\Transp& \zerobold
\end{bmatrix}
\begin{bmatrix}
\ubold \\
\lambdabold \\
\ubold_{I}
\end{bmatrix}
=
\begin{bmatrix}
\fbold\\
\zerobold\\
\zerobold
\end{bmatrix}
\end{equation}
by defining the following block matrices and vectors:
\begin{align}
\Bbold = \diag \begin{bmatrix} \Bbold_{\star}^{(1)} \dots \Bbold_{\star}^{(n_{p})} \end{bmatrix} \\
\Lbold = \begin{bmatrix} \Lbold_{\star}^{(1)} \\ \vdots \\
\Lbold_{\star}^{(n_{p})} \end{bmatrix}, \quad \lambdabold =
\begin{Bmatrix} \lambdabold_{\star}^{(1)} \\ \vdots \\
\lambdabold_{\star}^{(n_{p})} \end{Bmatrix}
\end{align}
with $\star = S,F$ depending on the physics associated to
substructure $p$, i.e., (S) for a structure modeled using the FEM or (F) for an acoustic fluid approximated with the BEM.

Finally, we are interested in solving the problem first for the interface.
This can be done obtaining the subdomain displacements from the
first row of (\ref{LLM_system}):
\begin{equation}
\ubold = \Kbold^{-1}(\fbold - \Bbold \lambdabold)
\end{equation}
 and substituting into the second row to obtain the following flexibility system:
\begin{equation} \label{eqn_ABETI}
\begin{bmatrix}
\Fbold_{bb}   & \Lbold \\
\Lbold\Transp& \zerobold
\end{bmatrix}
\begin{bmatrix}
\lambdabold\\
\ubold_{I}
\end{bmatrix}
=
\begin{bmatrix}
\bbold\\
\zerobold
\end{bmatrix}
\end{equation}
with $\Fbold_{bb} = \Bbold\Transp\Kbold^{-1}\Bbold$ and
$\bbold=\Bbold^T\Kbold^{-1}\fbold$.

\subsection{Test: Mortar-MLL comparisson}\label{section33}
\begin{figure}
\begin{center}
\includegraphics[width=7cm]{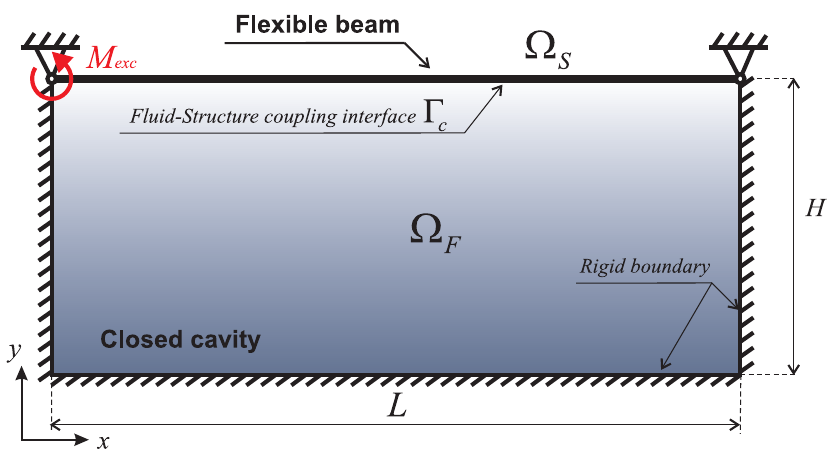} \\
\textbf{\small (a)}\\
\text{ }\\
\includegraphics[width=7.5cm]{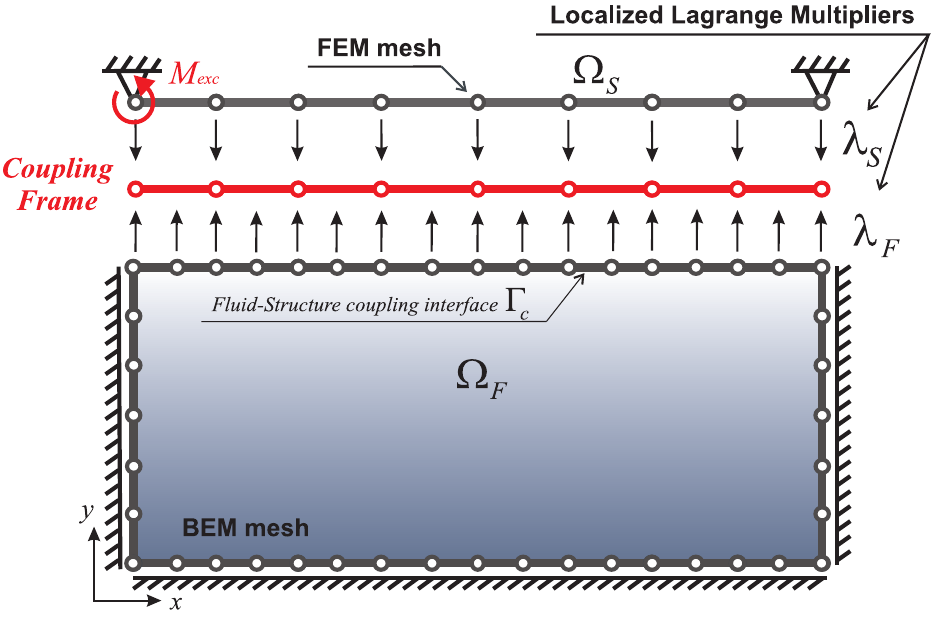}\\
\textbf{\small (b)}\\
\caption{Acoustic cavity with a flexible wall and harmonic
excitation. Problem definition (a) and BEM-FEM subdomains
coupled using LLM (b).} \label{ej1_modelo}
\end{center}
\end{figure}

The coupling possibilities of Mortar and LLM methodologies are
studied and compared by solving the following test problem taken from \cite{Sandberg08b}: a two
dimensional $L\times H$ cavity ($L=10$ $m$ and $H=4$ $m$) with
one flexible side (see Figure  \ref{ej1_modelo}). The flexible wall is a beam that is simply supported on both edges of the cavity and is modeled using Euler-Bernoulli beam elements. The properties of this structural domain are: Young module $E=$ $2.1\times10^{11}$ $Pa$, section inertia
$I=$ $1.59\times10^{-4}$ $m^4$, cross section area $A=0.02$ $m^2$ and a mass per unit length $m_s=50$ $kg/m$. The remaining three sides of the cavity are reverberant walls where homogeneous
Neumann boundary conditions are applied ($v_n=0$). The fluid is
water with $c_F=1500m/s$ and $\rho_F=1000kg/m^3$. The sketch of this problem is
presented in Figure \ref{ej1_modelo}(a) where a harmonic bending moment $M_{exc}=$ $M_{o}e^{i\omega t}$ is applied in one edge. In Figure \ref{ej1_modelo}(b) an scheme of the meshes and the coupled BEM-FEM subdomains using LLM is presented.

\begin{figure}
\begin{center}
\includegraphics[height=7cm,angle=0]{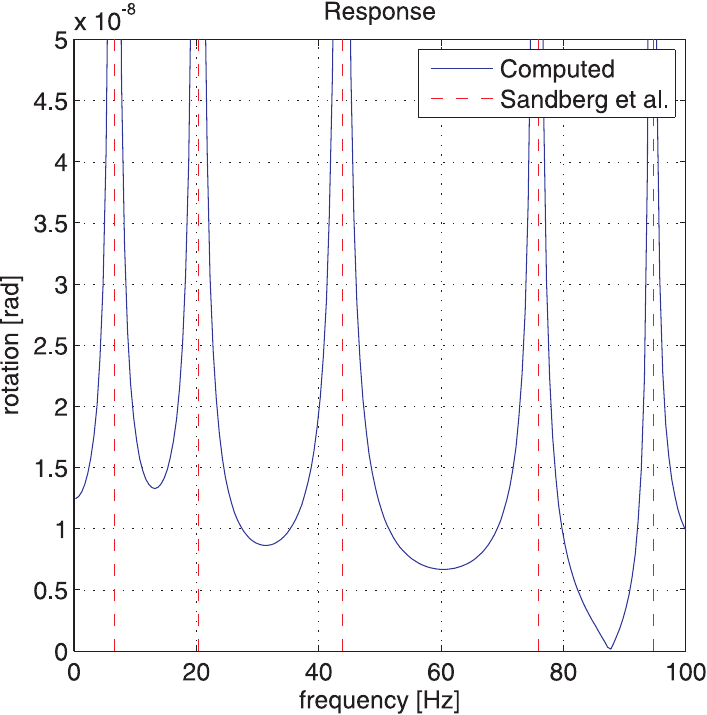}
\caption{Transfer function of the cavity problem for the beam rotation at $x=L$. Natural frequencies computed by Sandberg et al. \cite{Sandberg08b} using a FEM-FEM coupling method.}\label{ej1_repuestar}
\end{center}
\end{figure}

\begin{figure*}
\begin{center}
\includegraphics[width=6.5cm]{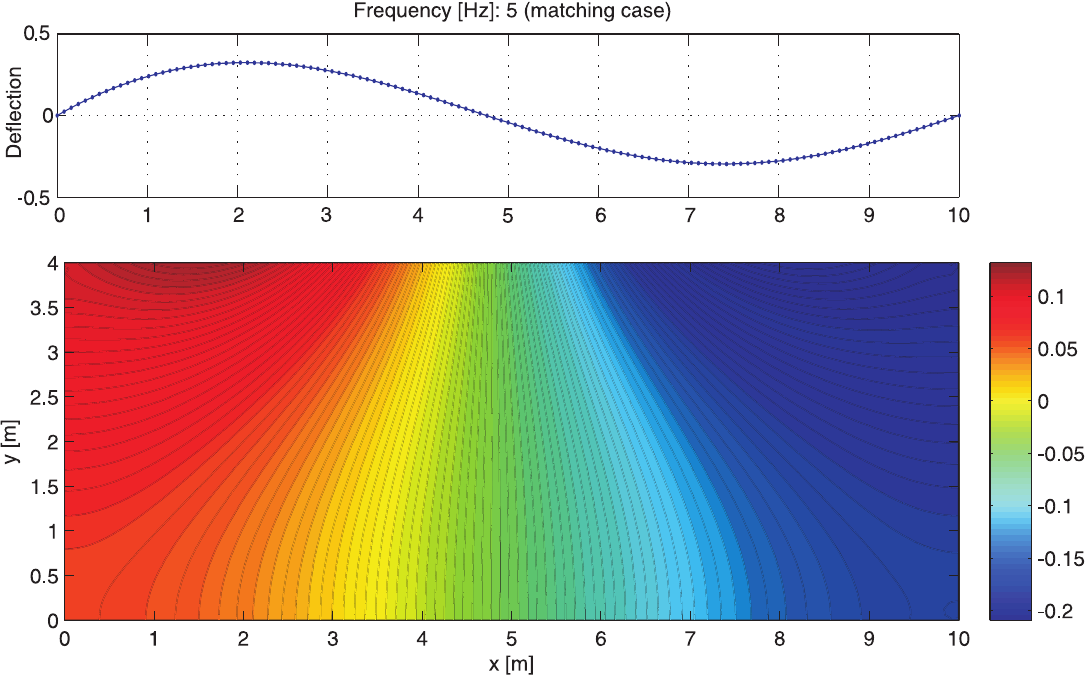}
\includegraphics[width=6.5cm]{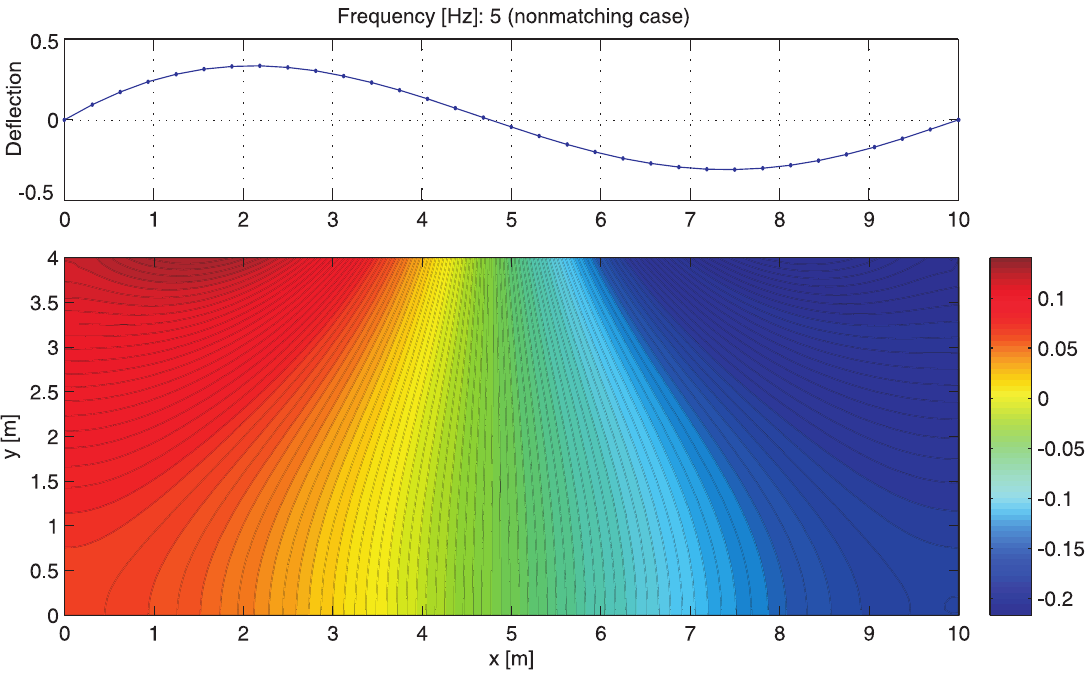}\\
\textbf{\small (a)}\\
\text{ }\\
\includegraphics[width=6.5cm]{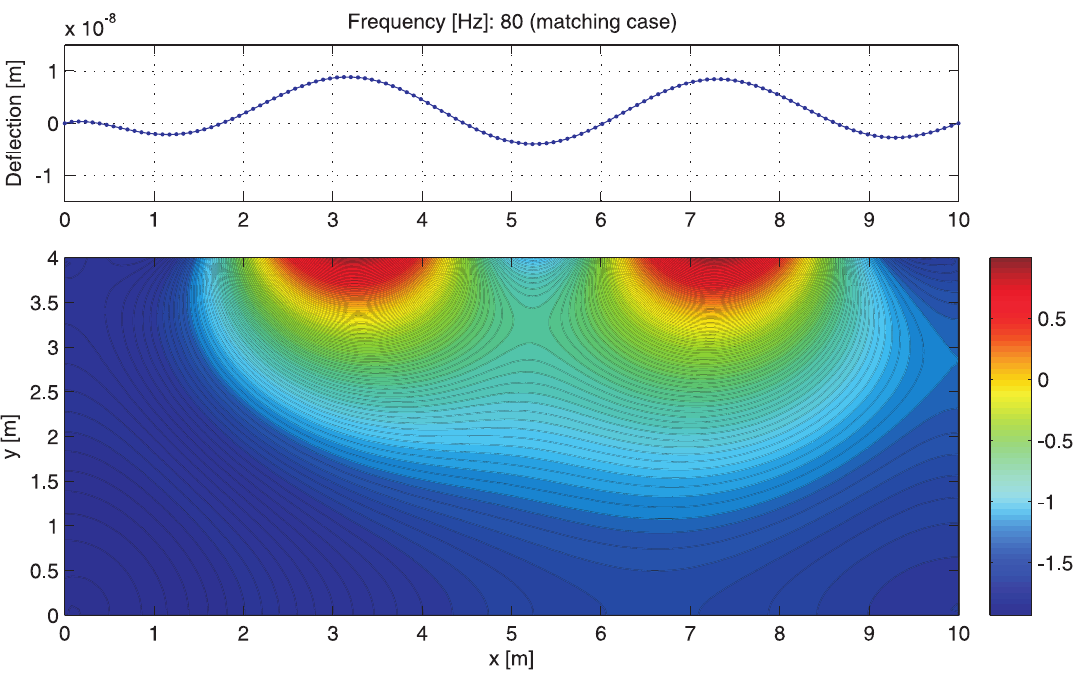}
\includegraphics[width=6.5cm]{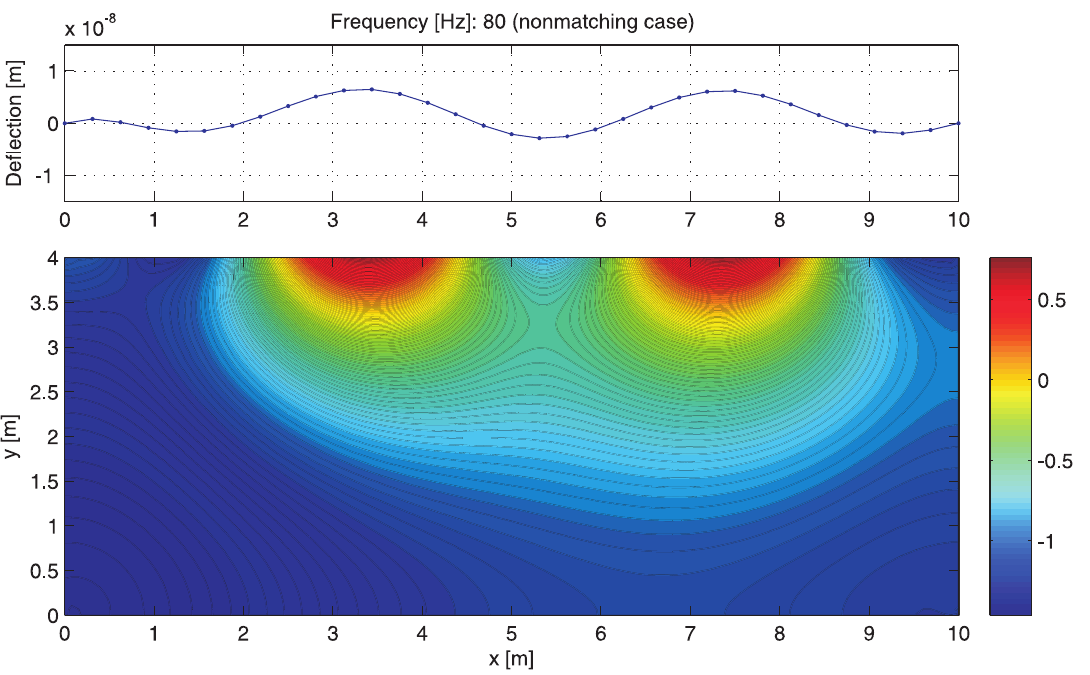}\\
\textbf{\small (b)} \caption{Acoustic cavity. Deflection of the flexible wall and
fluid pressure field for different excitation frequencies. Results
for $5$ Hz (a) and $80$ Hz (b) using matching meshes (left) and
non-matching meshes (right).}\label{ej1_SoluDomain}
\end{center}
\end{figure*}

Figure \ref{ej1_repuestar} shows the beam rotation at $x=L$ as a
function of the excitation frequency. The results coincide with the natural
frequencies obtained by Sandberg et al. in \cite{Sandberg08b}.
In Figure \ref{ej1_SoluDomain} it can be observed the flexible
wall deflection and fluid pressure due to excitations of $5$ Hz and $80$ Hz, for matching Figure \ref{ej1_SoluDomain}(a) and non-matching Figure \ref{ej1_SoluDomain}(b) meshes.

\begin{figure*}
\begin{center}
\includegraphics[width=12cm]{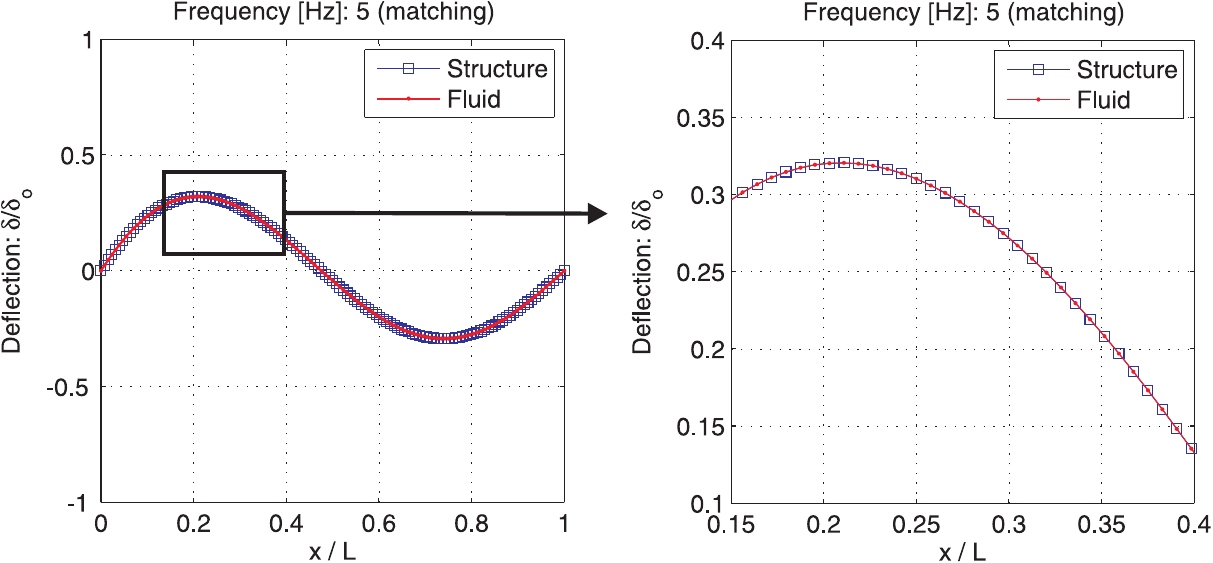} \\
\hspace{0.4cm}\textbf{\small (a)}\hspace{6cm} \textbf{\small (b)}\\
\text{ } \\
\includegraphics[width=12cm]{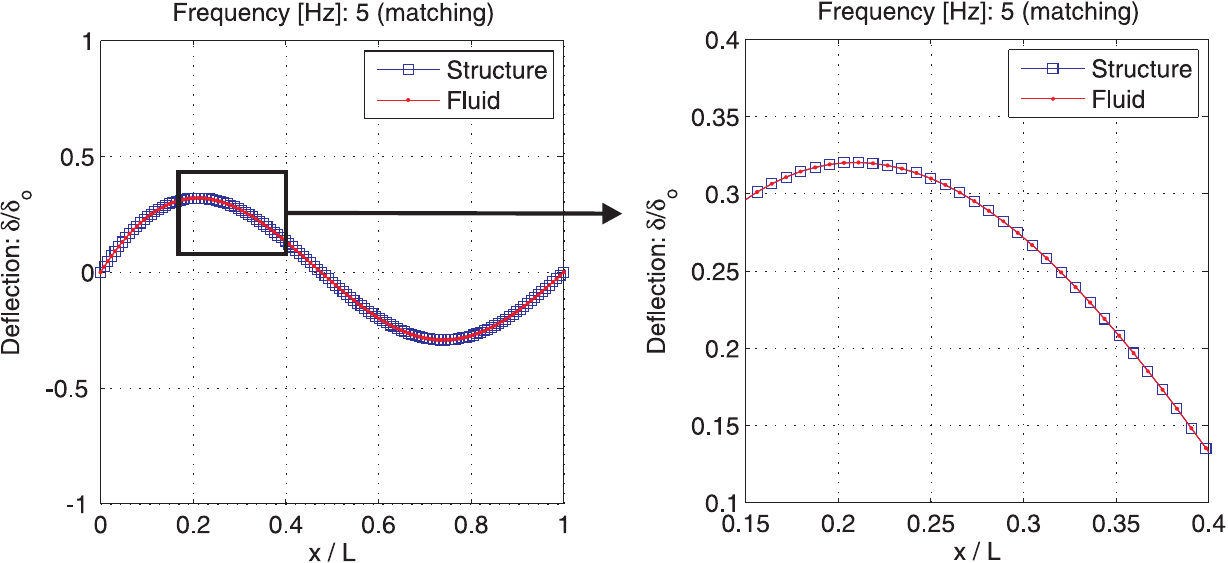} \\
\hspace{0.4cm}\textbf{\small (c)} \hspace{6cm} \textbf{\small (d)}\\
\caption{Beam deflection and interface fluid displacements for a
bending moment with frequency $\omega=5$ Hz. Interface coupling using matching meshes
with mortar method (top) and LLM (bottom).}\label{ej1_match_MortarLLM}
\end{center}
\end{figure*}

\begin{figure*}
\begin{center}
\includegraphics[width=12cm]{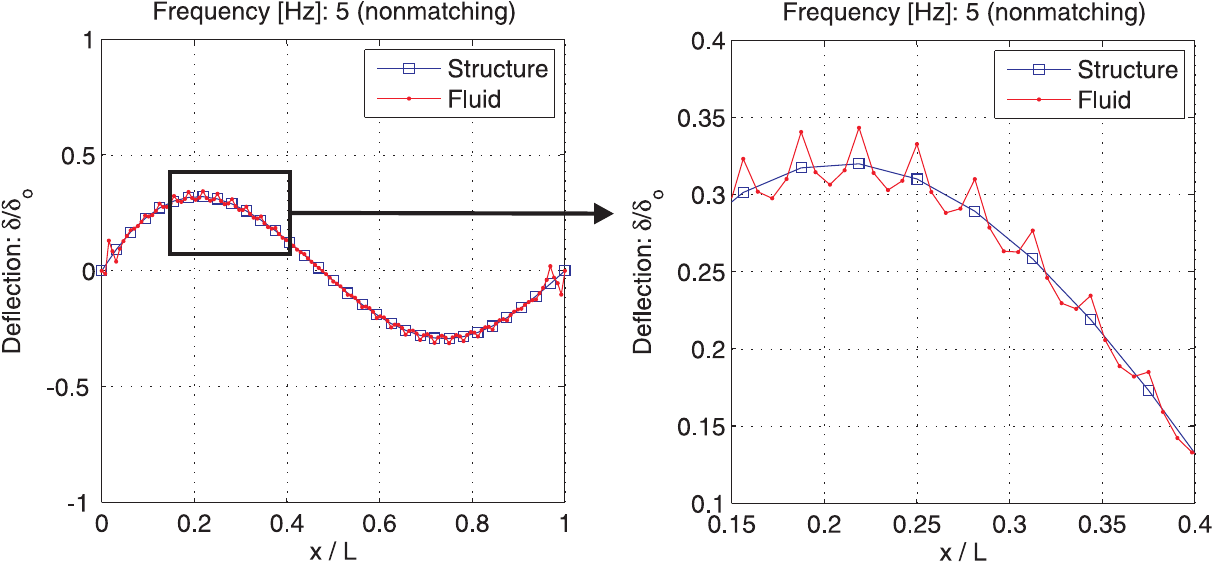} \\
\hspace{0.4cm}\textbf{\small (a)}\hspace{6cm} \textbf{\small (b)}\\
\text{ } \\
\includegraphics[width=12cm]{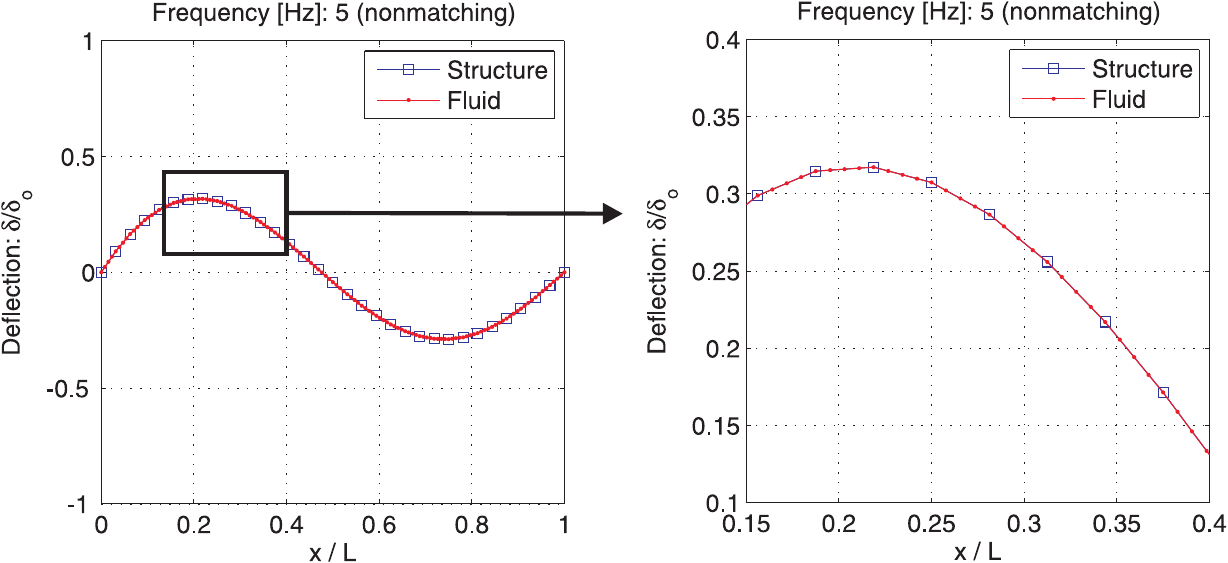} \\
\hspace{0.4cm}\textbf{\small (c)} \hspace{6cm} \textbf{\small (d)}\\
\caption{Beam deflection and interface fluid displacements for a
bending moment of frequency $\omega=5$ Hz. Interface coupling using non-matching meshes with mortar method (top) and LLM (bottom). In the mortar case, artifacts appear as a consequence of imposing the displacement compatibility condition
in a weak sense (top-right).}\label{ej1_nonmatch_MortarLLM}
\end{center}
\end{figure*}

The coupling interface displacements of the structure and the
fluid obtained using matching meshes with the mortar method and the LLM method are presented in
Figure \ref{ej1_match_MortarLLM}, for an harmonic excitations of $5$ Hz. As it is observed in
Figure \ref{ej1_match_MortarLLM}(a-b) and Figure \ref{ej1_match_MortarLLM}(c-d)
both methodologies present the same coupling behavior using
matching meshes at the interfaces. However, Figure \ref{ej1_nonmatch_MortarLLM} presents the coupling interface displacements computed using nonmatching meshes with the mortar and the LLM method. It can be observed in Figure \ref{ej1_nonmatch_MortarLLM}(a-b) that we obtain innacuracies in the form of wriggles on the fluid displacements using the mortar method.
The appearance of these artifacts in the mortar solution is attributed to the use of different displacement interpolations for the fluid (linear shape functions) and the structure (Hermite polynomials). One main characteristic of mortar methods is that the condition of pointwise continuity across the interface is replaced by a weak form and this standard primal approach is suboptimal when mixed finite element discretizations are used \cite{Wohlmuth00}.
For this reason, LLM method is going to be used in the nsBE-FETI
methodology.

\section{Iterative solution strategy for the interface problem}\label{sec_5}
The solution strategy presented to solve the flexibility system
obtained for the FSI localized Lagrange multipliers formulation \eqref{eqn_ABETI} uses a projection of the interface solution vector in the form:
\begin{equation} \label{multiplierDecomp}
\lambdabold = \mathcal{P}\lambdabold_d
\end{equation}
with the symmetric projector defined as:
\begin{equation}
\mathcal{P} = \mathbf{I} - \Lbold (\Lbold^T\Lbold)^{-1} \Lbold^T
\end{equation}
satisfying the condition: $\mathcal{P} \Lbold = \zerobold$.

Substitution of this decomposition into the flexibility formulation
of the interface problem \eqref{eqn_ABETI} yields the following equation:
\begin{equation}
\mathcal{P}\Fbold_{bb}\mathcal{P} \lambdabold_d= \mathcal{P}\bbold
\end{equation}
So the projected residual is then finally given by:
\begin{equation}
\rbold = \mathcal{P} (\bbold - \Fbold_{bb} \mathcal{P}
\lambdabold_{d}) \label{eqn_residBETI}
\end{equation}
equation that is solved for $\mathcal{P} \lambdabold_{d}$.

Because the non-symmetrical BEM-FEM interface problem is usually very large in practical applications, Krylov's iterative schemes for non-symetrical systems like Bi-CGSTAB and GMRES are prefered for the minimization of residual \eqref{eqn_residBETI}. The authors introduced in \cite{Gonzalez12c} a projected Bi-CGSTAB algorithm for non-symmetrical BEM problems in statics. This projected Bi-CGSTAB iterative scheme is generalized to dynamics in Table\ref{CBiCGStab} for the proposed nsBE-FETI formulation.

The proposed preconditioners for the fluid and the structure are
extensions of the well-known lumped and Dirichlet preconditioners of the standard FETI and AFETI algorithms. These preconditioners are calculated in a domain-by-domain basis as:
\begin{equation}
\tilde{\Fbold}_{bb}^{+} =
\begin{cases}
\barKbold_{bb} & \text{(FEM subdomain)} \\
\rho_F\omega^2\Mbold_{bb}\Hbold_{bb}^{-1}\Gbold_{bb} &\text{(BEM
subdomain)}
\end{cases}
\label{BEMpreconditioner}
\end{equation}
where subscript $(bb)$ refers to boundary extraction, i.e. pre and
post multiplication by $\Bbold\Transp $ and $\Bbold$ respectively.

Before an iterative method can be used to solve equation
\eqref{eqn_LLMsystem}, a scaling of the variables based on
\cite{Li11} should be applied to improve the condition number of the system.
Denoting $\Lambda=E\nu/((1+\nu)(1-2\nu))$, the scaled
displacements are
$\tilde{\dbold}_s=\sqrt{\rho\omega^{2}\Lambda}~\ubold_{S}$ and
$\tilde{\dbold}_f=\sqrt{\rho\omega^{2}\Lambda}~\ubold_{F}$, and
the scaled Lagrange multipliers $\tilde{\lambdabold}=1/\sqrt{\rho\omega^{2}\Lambda}~\lambdabold$.

\begin{table}
\centering
\begin{tabular}{|l|}
\hline
(I) Initialize: \\
\quad $\lambdabold_0$, $\rbold_0  = \mathcal{P} \left( {\bbold   - \Fbold_{bb} \lambdabold_0}\right)$ \\
\quad $\xbold_0  = \zerobold$, \hspace{0.2cm} $\pbold_0  = {\zerobold}$ \\
(II) Iterate $i=1,2,3...$ until convergence: \\
\quad Compute: \\
\qquad $\pbold_i  = \rbold_{i-1}  + \omega _i (\pbold_{i-1}  - \alpha _{i-1} \xbold_{i - 1})$ \\
\qquad with  $\pbold_{1}= \rbold_{0}$,  $\beta _i  = (\mathbf{\hat{r}}_0^{*}\rbold_{i - 1})$ \\
\qquad and $\omega _i  = {{\beta _i \gamma _{i - 1} }/({\alpha _{i - 1} \beta _{i - 1} }})$ \\
\quad Precondition: \\
\qquad $\mathbf{a}_i = {\bf{\tilde F}}^ +_{bb} \pbold_i$ \\
\quad Projection: \\
\qquad $\mathbf{z}_i = \mathcal{P} \mathbf{a}_i$ \\
\quad Compute: \\
\qquad $\ubold_i = \rbold_{i - 1}  - \gamma _i\xbold_i$, \\
\qquad with $\bbold_i = \Fbold_{bb}\mathbf{z}_i$, $\xbold_i  =\mathcal{P} \bbold_i$ and $\gamma _i  = {{\beta _i }/({\mathbf{\hat{x}}_i^{*}\rbold_0}})$ \\
\quad Precondition:\\
\qquad $\mathbf{c}_i = {\tilde{\Fbold}}^ +_{bb} \ubold_i$ \\
\quad Projection: \\
\qquad $\ybold_i = \mathcal{P} \mathbf{c}_i$ \\
\quad Update solution: \\
\qquad $\lambdabold_{i}=\lambdabold_{i-1}+\gamma _{i}\mathbf{z}_{i}+\alpha_{i}\ybold_{i}$ \\
\qquad with $\Gbold_i = \Fbold_{bb} \ybold_i$, $\mathbf{w}_i  = \mathcal{P} \Gbold_i$ \\
\qquad and $\alpha _i  = ({\mathbf{\hat{w}}_i^{*}\ubold_i })/({{\mathbf{\hat{w}}_i^{*}\mathbf{w}_i}})$ \\
\quad Update residual: \\
\qquad $\rbold_{i} = \ubold_{i}  - \alpha _i \mathbf{w}_i$ \\
(III) If ${{\left\| {\rbold_i } \right\|}/{\left\| {\rbold_0 } \right\|}}>\epsilon $ , $i\leftarrow i+1$ return to step (II) \\
\\ \hline
\end{tabular}
\caption{Complex BiCGStab algorithm with preconditioning and
projection used to minimize the residual of equation
\eqref{eqn_residBETI}.\label{CBiCGStab}}
\end{table}

\section{Numerical Results}
The possibilities of the proposed methodology are demonstrated in
this section, where three representative examples are investigated
solving the flexibility equation \eqref{eqn_ABETI} using the nsBE-FETI
iterative algorithm. The influence of different factors in the
convergence of the nsBE-FETI algorithm like the
number of elements per subdomain, frequency of the excitation
and presence of non-matching interfaces, are examined.

\subsection{Acoustic cavity with a flexible wall}
This first example revisits the problem presented in section
\ref{section33} (see Figure \ref{ej1_modelo}). A series of cases
using BEM-FEM matching meshes are first solved. Fluid domain is
discretized using linear boundary-element meshes with $L/h=32$,
$64$, $128$, and $256$ divisions at the interface and the
structure is discretized using two-node Euler-Bernoulli beam
elements. Two frequencies excitation of $5$ $Hz$ and $80$ $Hz$ are
considered to study the influence of the frequency in the convergence.

\begin{figure}
\begin{center}
\includegraphics[height=6cm,angle=0]{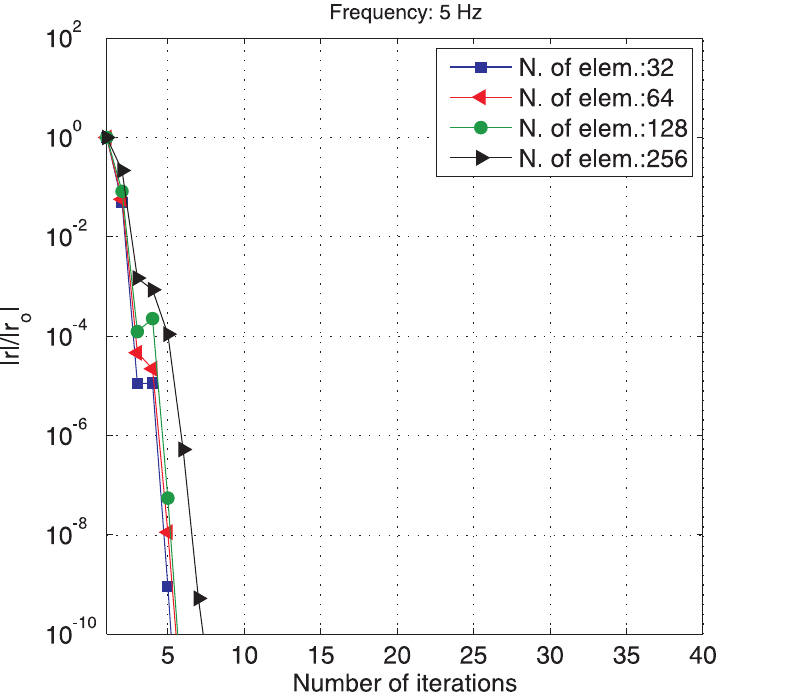}\\
\textbf{\small (a)}\\
\text{ }\\
\includegraphics[height=6cm,angle=0]{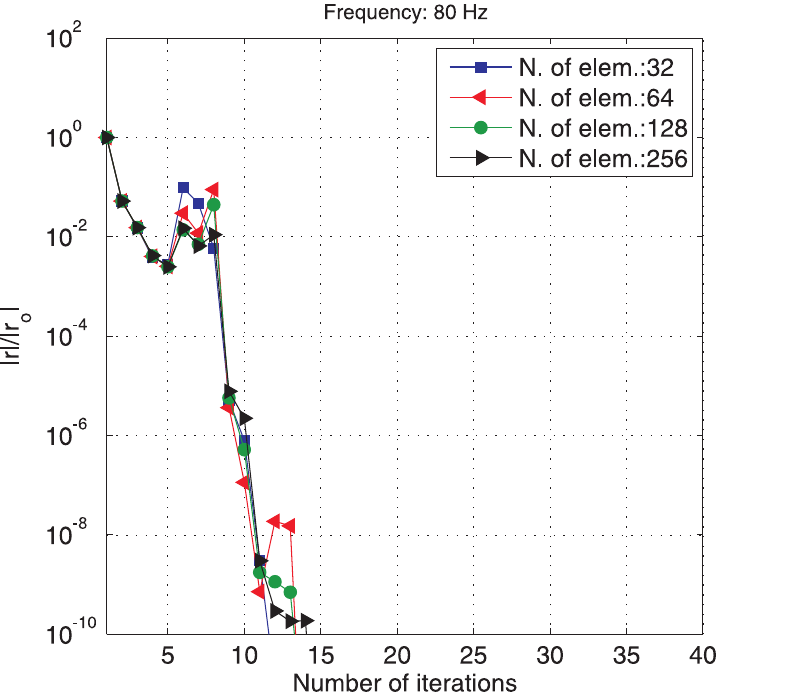}\\
\textbf{\small (b)}\\
\caption{BiCGSTAB error evolution for: $5$ Hz (a) and $80$ Hz (b),
considering a LLM coupling of matching meshes.}
\label{ej1_LLM_matching_fig}
\end{center}
\end{figure}

\begin{figure}
\begin{center}
\includegraphics[height=6cm,angle=0]{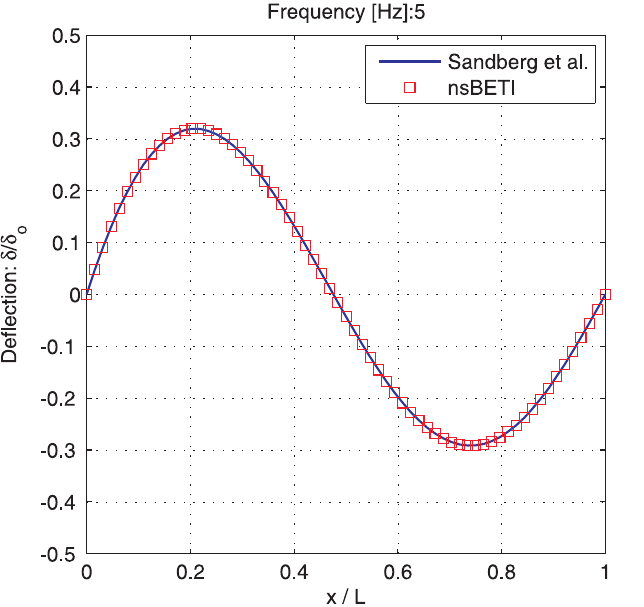}\\
\textbf{\small (a)}\\
\text{ }\\
\includegraphics[height=6cm,angle=0]{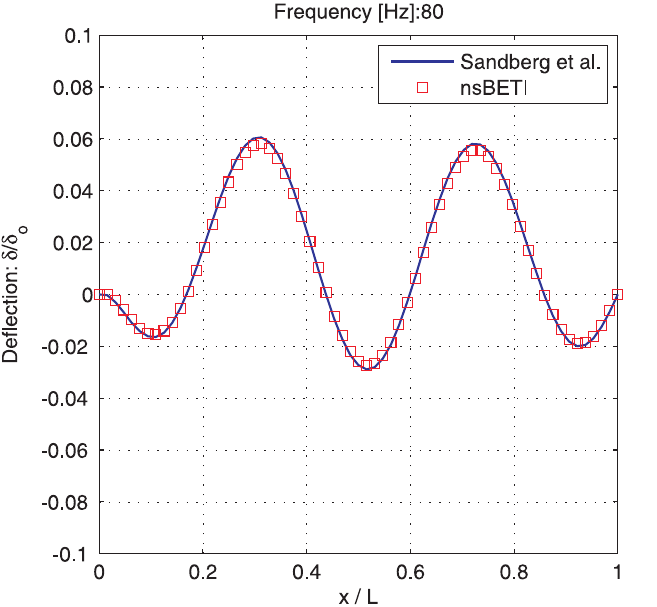}\\
\textbf{\small (b)}\\
\caption{Beam deflectiosn due to harmonic excitation of different frequencies: $5$ Hz (a) and $80$ Hz (b), considering a LLM coupling of matching
meshes.}\label{ej1_LLM_matching_fig2}
\end{center}
\end{figure}

Table \ref{ej1_tabla1} shows the number of iterations needed by
the projected Bi-CGSTAB algorithm to solve these coupled problems with a tolerance of $10^{-10}$. Figures \ref{ej1_LLM_matching_fig}(a) and
\ref{ej1_LLM_matching_fig}(b) show the convergence evolution for a low  excitation frequency of $5$ Hz and a higher frequency of $80$ Hz with the number of iterations needed by the algorithm to solve these problems. For the cases considered, it can be observed that an exponential increase of the type $L/h = 2^n$ translates into a constant number of iterations for both excitation frequencies. The
difference in the iterations number between $5$ Hz
and $80$ Hz is due to the complexity of the solution, as Figure \ref{ej1_LLM_matching_fig2} shows.

\begin{figure}
\begin{center}
\includegraphics[height=6cm,angle=0]{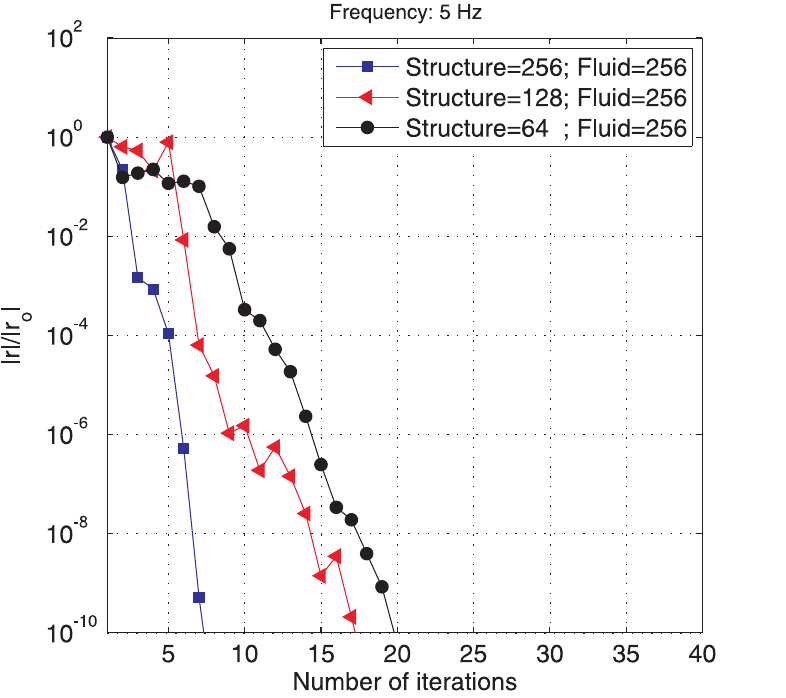}\\
\textbf{\small (a)}\\
\text{ }\\
\includegraphics[height=6cm,angle=0]{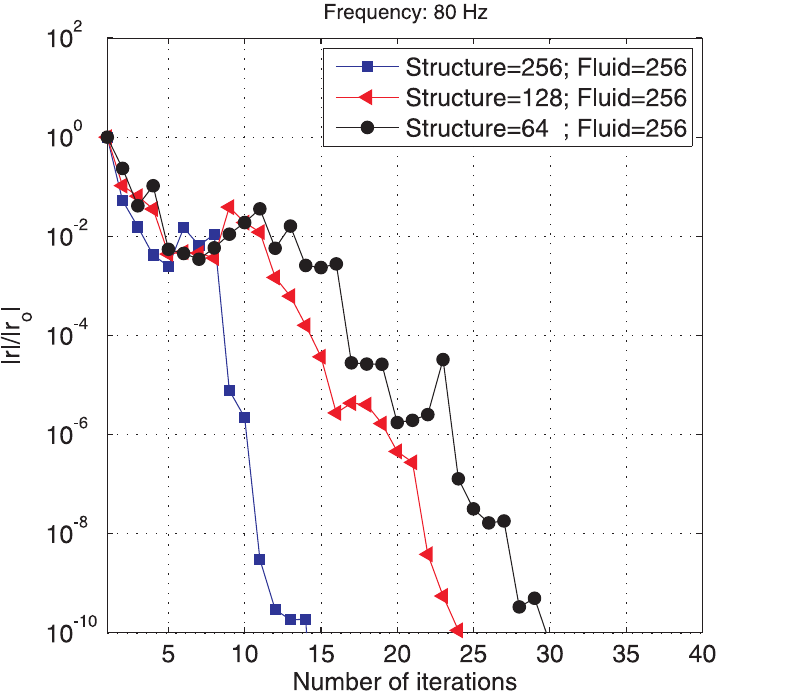}\\
\textbf{\small (b)}\\
\caption{Bi-CGSTAB error evolution considering non-matching meshes and harmonic excitations of: $5$ Hz (a) and $80$ Hz (b).}%
\label{ej1_LLM_nonmatching_fig}
\end{center}
\end{figure}

Finally, the non-matching case is considered changing the
discretization of the structure to produce dissimilar meshes at
the interface. Figure \ref{ej1_LLM_nonmatching_fig} presents error evolutions for $5$ Hz (Figure \ref{ej1_LLM_nonmatching_fig}(a)) and $80$ Hz (Figure \ref{ej1_LLM_nonmatching_fig}(b)). The results are obtained for structural meshes ranging from $L/h=64$ (highly
non-matching case) to $256$ (matching case) maintaining the mesh of the fluid fixed with $L/h=256$ divisions. It
is noted that the introduction of dissimilar meshes, maintaining a
constant $(L/h)_{max}$, slightly increases the
number-of-iterations needed by nsBE-FETI to solve the problem for low
and high excitation frequencies. The experiment is repeated for
the BEM-FEM case, see Table \ref{ej1_tabla2}, presenting similar
results. As a summary, in the matching case, the
convergence of nsBE-FETI is governed by $(L/h)_{max}$, but
the introduction of non-matching interfaces destructs this
property producing a negative effect in the convergence that is
controlled by the interface mesh-dissimilarity parameter
$h_{max}/h_{min}$. However, for the cases studied, the impact of a
non-matching interface is limited and does not significantly
affect the algorithm convergence.

\begin{table}
\begin{center} \begin{tabular}{ccccc}
\hline
$L/h$  & $n_{e}$  & $n_{e}$ & Iterations & Iterations \\
& BEM & FEM & $5Hz$ & $80Hz$ \\
\hline
32      & 96           & 32   &  6   &  12   \\
64      & 192          & 64   &  6   &  14   \\
128     & 384          & 126  &  6   &  14   \\
256     & 768          & 256  &  8   &  15   \\
\hline
\end{tabular} \end{center}
\caption{Acoustic cavity problem with BEM-FEM matching interface. Number of iterations for a constant normalized residual with different mesh sizes.} \label{ej1_tabla1}
\end{table}

\begin{table}
  \centering
    \begin{tabular}{cccc}
    \hline
$L/h$ & $L/h$ & Iterations & Iterations \\
FEM & BEM  & $5 Hz$ & $80 Hz$ \\
    \hline
        64      & 256      &  8     &  15    \\
        128     & 256      &  18    &  25    \\
        256     & 256      &  20    &  30    \\
    \hline
    \end{tabular}
  \caption{Acoustic cavity problem with BEM-FEM non-matching interface. Number of iterations for a normalized residual of $10^{-10}$ with different mesh sizes.}\label{ej1_tabla2}
\end{table}

\subsection{Rectangular duct with closed outlet}
\begin{figure*}
\begin{center}
\includegraphics[width=10cm]{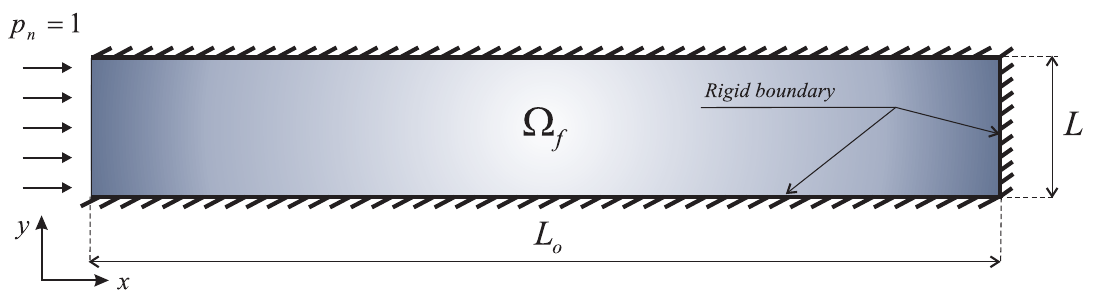}\\
\textbf{\small (a)}\\
\text{ }\\
\includegraphics[height=4cm,angle=0]{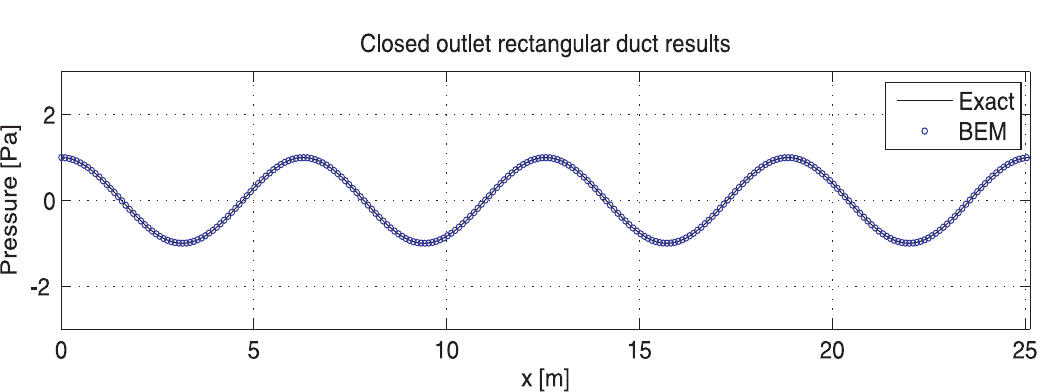}\\
\includegraphics[height=4cm,angle=0]{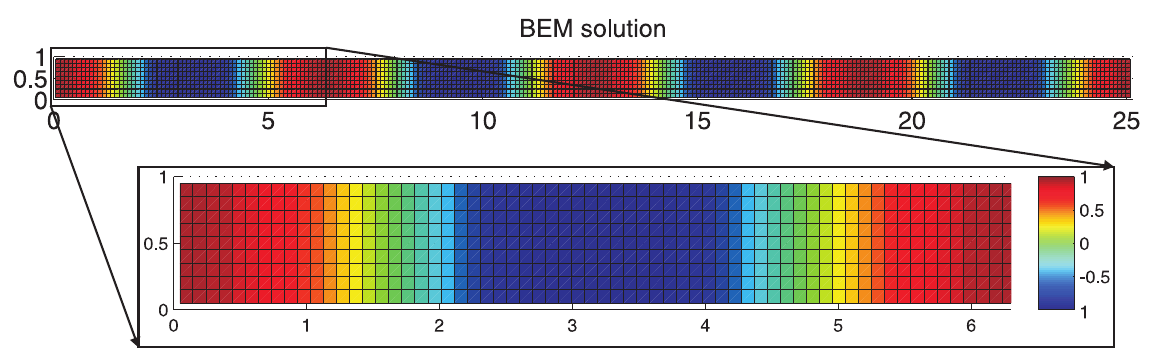}\\
\textbf{\small (b)}\\
\caption{Rectangular duct with closed outlet. (a) Problem description, dimensions and boundary conditions. (b) Distribution of the
fluid pressure in the longitudinal direction compared with the
analytical solution.} \label{ej2_examp_solucion2}
\end{center}
\end{figure*}

Next example considers a simple rectangular duct with a closed
outlet as represented in Figure \ref{ej2_examp_solucion2}(a). The closed outlet is located at $x=L_o$ and assumed to be a rigid wall from $y=0$ to $y=L$, and the inlet has a complex pressure ($p=$ $p_{o}e^{i\omega t}$) prescribed at $x=0$. The fluid is water as in the previous example. The wave number is set to $k=1$, the length of the duct in the $x$-direction is $L_o=$ $8\pi$ $m$ and the section height is $L=$ $1m$. Figure \ref{ej2_examp_solucion2}(b) presents the solution in terms of resulting pressure distribution on the field points.

\begin{figure*}
\begin{center}
\includegraphics[width=14cm]{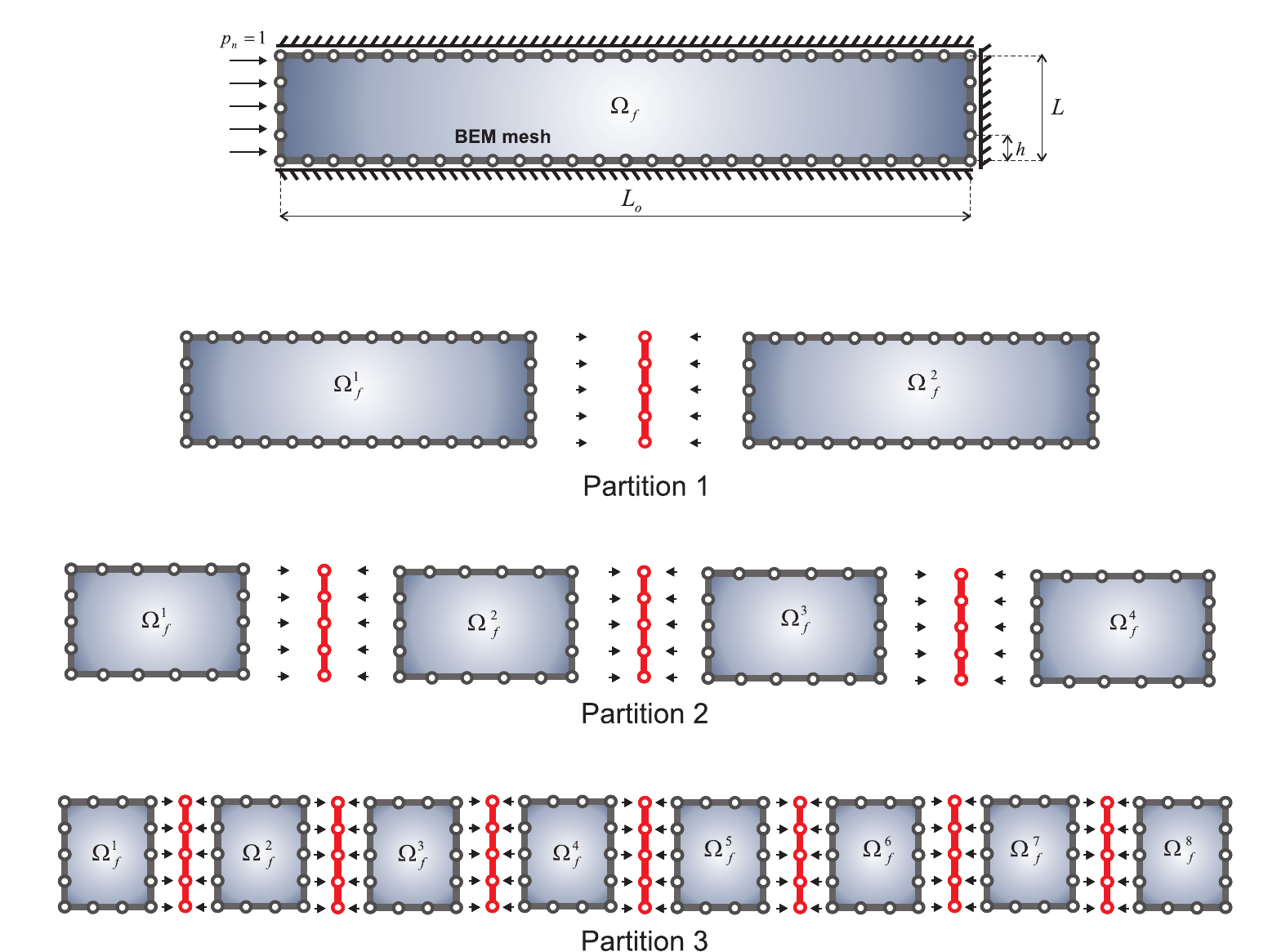}
\caption{Rectangular duct with closed outlet. Partitioning of the
fluid domain into $n_{p}$=2, 4 and 8 subdomains connected with
localized Lagrange multipliers.} \label{ej2_partition}
\end{center}
\end{figure*}

\begin{figure}
\begin{center}
\includegraphics[height=6cm,angle=0]{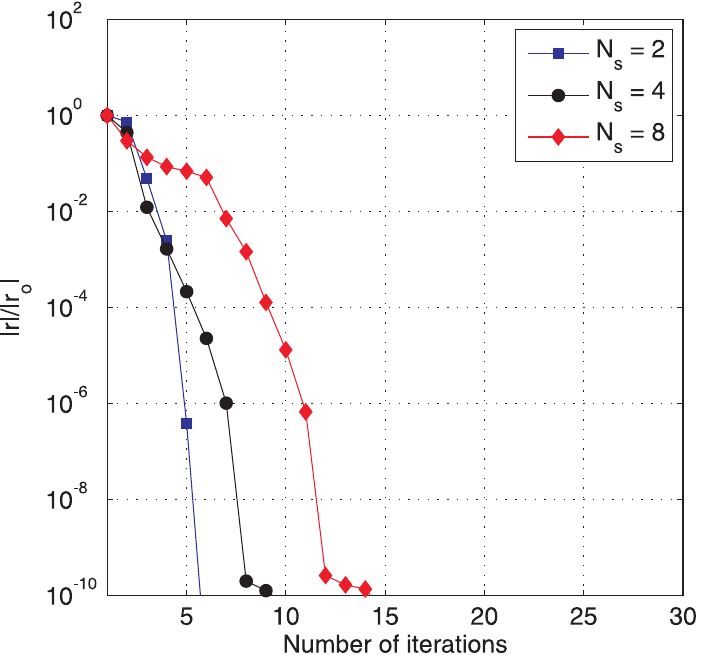}
\caption{Duct problem. Evolution of the residual for different number of partitions.}%
\label{ej2_error}
\end{center}
\end{figure}

The duct is partitioned transversally into $N_s=$ $2$, $4$ and $8$
subdomains, discretized using linear boundary elements of fixed size $L/h=$ $10$, maintaining matching interfaces (Figure \ref{ej2_partition}). The objective of this test is to demonstrate that, maintaining the element size, the convergence of nsBE-FETI is not considerably affected by the number of partitions. Table \ref{ej2a_tabla1} contains a summary of convergence results, and Figure \ref{ej2_error} shows the convergence history of this particular case. It is observed a small effect of $N_s$ in the convergence.

In the next experiment, the number of subdomains is fixed to $N_s=$ $4$ and the problem solved for different discretizations with $L/h=$ $10$, $20$ and $40$, using a total of $148$, $296$ and $592$
elements-per-subdomain. Convergence-rate results are
summarized in Table \ref{ej2a_tabla2}, demonstrating a logarithmic
correlation between the number of iterations for convergence and the mesh-size ($L/h$) in the range of discretization-sizes
studied.

\begin{table}[t]
  \centering
    \begin{tabular}{cccc}
    \hline
    $L/h$      & $n_{p}$        & $n_{e}$            & Iterations       \\
    \hline
    10 &  2 & 552   & 6  \\
    10 &  4 & 592   & 10 \\
    10 &  8 & 672   & 15  \\
    \hline
    \end{tabular}
  \caption{Influence of the number of partitions ($N_{s}$) for a fixed mesh discretization ($L/h$) of the duct problem.}
  \label{ej2a_tabla1}
\end{table}

\begin{table}
  \centering
    \begin{tabular}{cccc}
    \hline
    $L/h$      & $n_{p}$        & $n_{e}$            & Iterations       \\
    \hline
    10  & 4  & 592  & 10  \\
    20  & 4  & 1184 & 14  \\
    40  & 4  & 2368 & 18  \\
    \hline
    \end{tabular}
  \caption{The number of partitions is fixed ($N_{s}=4$) and the number of elements ($N_{el}$) increases.}
  \label{ej2a_tabla2}
\end{table}

\subsection{Open problem: Scattering object with a flexible wall}
\begin{figure*}
\begin{center}
\includegraphics[height=5.8cm,angle=0]{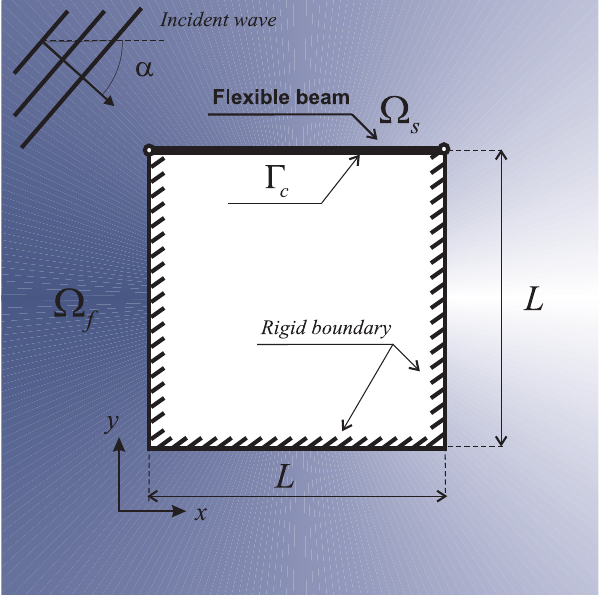}
\includegraphics[height=6cm,angle=0]{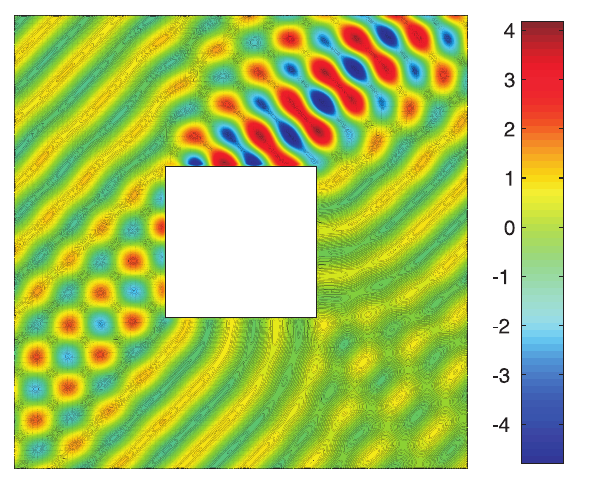} \\
\textbf{\small (a)} \hspace{5cm} \textbf{\small (b)}
\caption{Wave scattering produced by a monochromatic incident wave on a square obstacle with a flexible wall. (a) Problem definition. (b)
Real part of the total acoustic pressure for frequency $500Hz$ and incidence angle $\alpha=-\pi/4$.} \label{ej3_fig1}
\end{center}
\end{figure*}

Finally, in our last example we consider an open problem with a
square scattering object that has a flexible wall of length $L=$ $10$ $m$ (Figure \ref{ej3_fig1}(a)). The fluid is water and the structural domain presents the same properties than the first example. Our object is excited by a plane monochromatic wave of frequency $500$ $Hz$ and an incidence angle $\alpha=\pi/4$ $rad$.

\begin{figure}
\begin{center}
\includegraphics[height=6cm,angle=0]{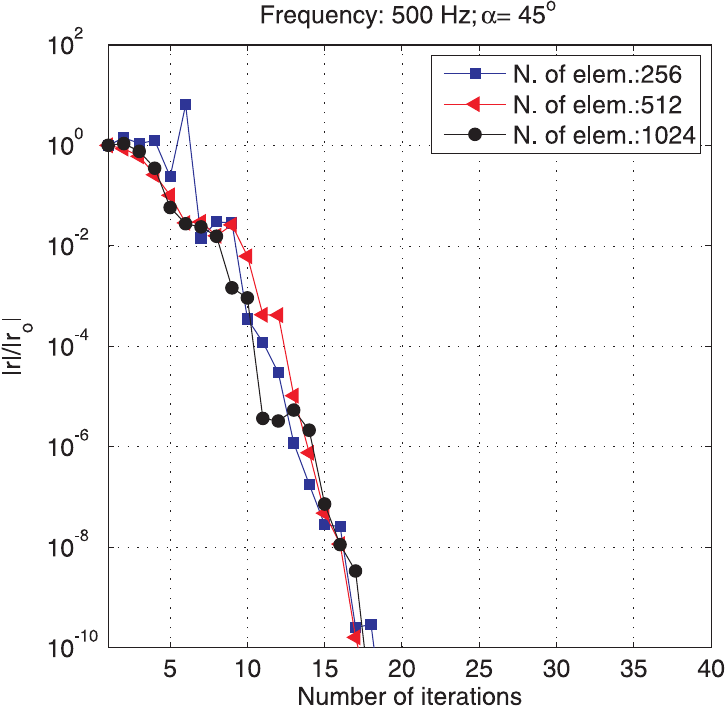}
\caption{Scattering problem. Convergence of the projected Bi-CGSTAB algorithm
for a frequency of $500$ Hz with different meshes.}
\label{ej3_solu}
\end{center}
\end{figure}

The domains are discretized using the same number of elements at
the coupling interface. Figure \ref{ej3_fig1}(b) shows the real
part of the total acoustic pressure around the object. Table \ref{ej3_tabla1} presents the
number of iterations needed by the projected Bi-CGSTAB algorithm to solve the coupled problem with a tolerance $\epsilon = 10^{-10}$ and Figure \ref{ej3_solu} shows the evolution of the residuals. It can be
observed a similar behaviour of the algorithm in this exterior problem than in the interior cases previously studied.

\begin{table}
\begin{center} \begin{tabular}{ccccc}
\hline
$L/h$       & $n_{e}$   & $n_{e}$       &  Iterations \\
& BEM & FEM & $500Hz$ \\
\hline
64      & 256          & 64   &  19 \\
128     & 512          & 128  &  18 \\
256     & 1024         & 256  &  18 \\
\hline
\end{tabular} \end{center}
\caption{Open problem with a BEM-FEM matching interface. Number of
iterations for a constant normalized residual with different mesh
sizes.}\label{ej3_tabla1}
\end{table}

\section{Summary and conclusions}\label{sec_6}
NsBE-FETI, a FETI-type formulation, has been extended to treat non-matching and non-symmetrical BEM-FEM acoustic FSI problems. This new formulation enjoys similar scalability properties than the classical FETI and symmetrical-BETI algorithms.

This resolution scheme is based on the LLM methodology which
allows to consider non-matching interfaces and preserves software
modularity. A comparison between LLM and the mortar scheme reveals
that the LLM method obtain a better interface displacements approximation than mortar for this kind of FSI problem: flexible wall
discretized using cubic beam elements coupled with an acoustics
fluid cavity, when highly dissimilar meshes are considered at the
interfaces.

Some scalability properties of the nsBE-FETI scheme have been studied considering different physics. First example was an interior acoustic problem with a flexible wall, where fluid and structure were discretized using matching and non-matching meshes. It was found that, in the matching case, convergence of nsBE-FETI algorithm is governed by the element size $(L/h)_{max}$ but the introduction of non-matching interfaces produces a negative effect in the convergence that is controlled by the interface mesh-dissimilarity parameter $h_{max}/h_{min}$. However, for the cases studied, the impact of a non-matching interface is limited and does not significantly affect the algorithm convergence for low frequencies.

In the second example, we modify the number of subdomains. Convergence studies reveal that for a fixed element size, the nsBE-FETI is not considerably affected by the number of partitions ($N_s$). Furthermore, when the number of subdomains is fixed and the problem solved for different discretizations, a logarithmic correlation between the number of iterations for convergence and the mesh-size used ($L/h$) is observed.

Finally, the last example presents an exterior FSI scattering problem where the same scalability behavior than in the interior ones could be observed. We can conclude that the proposed nsBE.FETI formulation equiped with the projected Bi-CGSTAB iterative solution algorithm presents good scalability properties for the solution of acoustic FSI problems.

\section{Acknowledgements}
This work was co-funded by the \emph{Ministerio de
Ciencia e Innovaci\'{o}n} (Spain), through the research projects DPI$2010$-$19331$, which is co-funded with the European Regional Development Fund (ERDF) (\emph{Fondo Europeo de Desarrollo Regional}, FEDER).


\clearpage

\bibliographystyle{unsrt}  



\end{document}